\documentclass{article}
\PassOptionsToPackage{numbers, compress}{natbib}

\usepackage[preprint]{style}

\usepackage{amsmath, amssymb}
\usepackage{amsthm}
\newtheorem{definition}{Definition}
\usepackage{dsfont}
\usepackage{algorithm}
\usepackage{booktabs}
\usepackage{algpseudocode}
\usepackage{graphicx}
\usepackage{mathtools}
\usepackage{multirow}
\usepackage{pifont}
\usepackage{subcaption}
\usepackage[dvipsnames]{xcolor}
\usepackage{tikz}
\usepackage{url}
\usepackage{hyperref}
\usepackage[misc]{ifsym}

\usepackage{tikz}
\usetikzlibrary{backgrounds, snakes, chains} 
\usetikzlibrary{arrows, decorations.markings, positioning}

\newcommand{\vx}{\boldsymbol x}
\newcommand{\vc}{\boldsymbol c}
\newcommand{\vb}{\boldsymbol b}
\newcommand{\va}{\boldsymbol a}
\newcommand{\vtheta}{\boldsymbol \theta}
\newcommand{\termdefn}[1]{\emph{#1}}

\definecolor{prettyblue}{rgb}{0.2, 0.2, 0.6}
\definecolor{cocoabrown}{rgb}{0.82, 0.41, 0.12}

\title{Machine Learning Augmented Branch and Bound for Mixed Integer Linear Programming
}

\newcommand{\email}[1]{\url{#1}}

\author{Lara Scavuzzo\thanks{Corresponding author} \\
        Delft University of Technology\\
        \email{L.V.ScavuzzoMontana@tudelft.nl} 
\And
        Karen Aardal \\
        Delft University of Technology\\
        \email{k.i.aardal@tudelft.nl}\\
\And
        Andrea Lodi \\
	Jacobs Technion-Cornell Institute\\
        Cornell Tech and Technion - IIT\\
	\email{andrea.lodi@cornell.edu} \\
\And
        Neil Yorke-Smith \\
	Delft University of Technology\\
	\email{n.yorke-smith@tudelft.nl}\\
}
\date{}

\begin{document}

\maketitle

\begin{abstract}
Mixed Integer Linear Programming (MILP) is a pillar of mathematical optimization that offers a powerful modeling language for a wide range of applications. 
During the past decades, enormous algorithmic progress has been made in solving MILPs, and many commercial
and academic software packages exist. Nevertheless, the availability of data, both from problem instances and from solvers, and the desire to solve new problems and larger (real-life) instances, trigger the need for continuing algorithmic development.
MILP solvers use branch and bound as their main component.  In recent years, there has been an explosive development
in the use of machine learning algorithms for enhancing all main tasks involved in the branch-and-bound algorithm, such as 
primal heuristics, branching, cutting planes, node selection and solver configuration decisions.
This paper presents a survey of such approaches, addressing the vision of integration of machine learning and mathematical optimization as complementary 
technologies, and how this integration can benefit MILP solving.
In particular, we give detailed attention to machine learning algorithms that automatically optimize some metric of 
branch-and-bound efficiency. 
We also address how to represent MILPs in the context of applying learning algorithms, MILP benchmarks and software.

\end{abstract}

\section{Introduction}
\label{sec:intro}
Mixed Integer Linear Programming (MILP) is a pillar of the field of mathematical optimization. Its widespread use covers multiple application domains such as transportation \cite{HandbookTransportation}, production planning \cite{PochetWolsey} and energy systems \cite{Wirtz2021}. The theory of modelling and solving MILPs is rich, and has fruitful intersections with several areas of mathematics and computer science \cite{NemhauserWolsey,Schrijver,Wolsey}. 
In practice, modern MILP solvers, which all use the tree search method {\em branch and bound (B\&B)} as the main component, are of very high quality and can solve, to provable optimality, problem instances that were considered as completely out of reach only a decade ago. 
This progress is to a large extent due to remarkable algorithmic advances in the past two decades (see, e.g., \citet{Achterberg2013} for an overview). 

Continuing these waves of algorithmic progress is key to unlocking new application domains and larger scales, as many areas of potential improvement still exist. For example, several algorithmic decisions that must be made within B\&B are made by heuristic rules that have been developed and tuned via computational studies to yield good average performance. For some of these decisions, we know that alternative heuristics exist that would lead to smaller search trees. These heuristics are, however, too time consuming in most cases. Recently, Machine Learning (ML) methodologies have been explored as a potential tool to mimic, or even improve, such time-consuming heuristics in order to gain the benefit of smaller search trees without the heavy computational burden. This strategy is particularly promising given the increase in data availability, not only in terms of problem instances but also data collected by the solver during the solution process. 



This survey addresses precisely the perspective of enhancing key components of B\&B by ML so as to 
integrate machine learning and mathematical optimization 
as complementary technologies, not as competing ones. We note that collecting data about the solving process and exploiting it to make informed decisions within the algorithm is standard practice in MILP solving. We survey methods that take this idea one step further, in the sense that the mapping between the collected data and the decision is not fixed a priori by an expert, but instead automatically constructed by `learning algorithms' that try to optimize a predefined metric of efficiency.

\paragraph{Scope.}
The survey by \citet{Bengio2021} presents the first developments of using machine learning in the context of optimization, with a focus on combinatorial optimization.
Our survey is zooming in on the B\&B framework in an MILP solver for solving 
a (class of) MILP.
The methods we discuss do not assume any particular problem structure a priori, in the sense that the structure is instead learned.
Our focus is on using ML for elements of B\&B where choices are being made, such as how to navigate the search tree, how to quickly find good feasible solutions, and how to improve linear relaxations.
We do not cover solution methods that replace B\&B, such as end-to-end ML approaches. By restricting the attention to the MILP context and its integration with ML, we are able to make a significant step forward in the characterization of choices like the problem representation and of specific aspects like benchmarks and software.


\paragraph{Readership.}
This survey is intended first of all for readers with some background in mathematical programming, and possibly with more limited background in machine learning.  We aim for it to be accessible to those interested in how ML can aid aspects of MILP solvers.  Highly detailed knowledge of how contemporary MILP solvers function is not required.  We assume some familiarity with B\&B.

\paragraph{Outline.}
The remainder of this introduction recalls the definition of an MILP, overviews the main components of a B\&B-based MILP solver and presents classical evaluation metrics that will be useful.  Section~\ref{sec:into_to_ml} introduces the main concepts of machine learning, with attention to those found most useful for B\&B.
Section~\ref{sec:tasks} proceeds to survey the tasks in MILP solving where ML can be useful.  Section~\ref{sec:representations} concerns representations of MILP, and in Section~\ref{sec:benchmarks} 
rounds out the survey by 
collating instance benchmarks and software used in the literature. Some final remarks and perspectives are given in Section~\ref{sec:conclusions}.

\subsection{Mixed Integer Linear Programming}
\label{subsec:intro_to_milp}

We follow the notation from Wolsey \cite{Wolsey}, to which the reader is referred to for a general introduction to MILP.  We are given a matrix $A\in \mathbb{Q}^{m\times n}$, vectors $\vc\in\mathbb{Q}^n$ and $\vb\in\mathbb{Q}^m$, and a partition $(\mathcal{A}, \mathcal{B}, \mathcal{C})$ of the variable index set $\{1,...,n\}$. A Mixed Integer Linear Program is the problem of finding
\begin{equation}
\label{eq:MILP}
    \begin{aligned}
 z^* = \min\ \  & \vc^T \vx  \\
 \text{subject to } & A \vx \geq \vb ,  \\
 & x_{j}\in \mathbb{Z}_{\geq 0} \hspace{3mm} \forall j\in \mathcal{
 A}, \\
 & x_{j}\in \{0,1\}\hspace{3mm} \forall j\in \mathcal{B}, \\
 & x_{j}\geq 0 \hspace{3mm} \forall j\in \mathcal{C}. \\
\end{aligned}
\end{equation}

We write $\mathcal{I}:=\mathcal{A}\cup \mathcal{B}$ to denote the set of integer and binary variables. An MILP is said to be mixed-binary if $\mathcal{A}=\emptyset$ and binary if $\mathcal{A}=\mathcal{C}=\emptyset$. Relaxing the integrality constraints for variables in $\mathcal{I}$ yields a Linear Program (LP), known as the LP relaxation of the MILP.

\subsection{MILP solvers}
\label{sec:solvers}
Modern exact solvers for MILP\footnote{Throughout the paper, we use the term `MILP solver' rather than the, sometimes more common, `MIP solver', the reason being that we especially focus on the learning developed for solving mixed integer \emph{linear} programs, not including the nonlinear case.  Nevertheless, most of the technology is used in the solvers for both linear and nonlinear programs.} implement the branch-and-bound algorithm. In our survey, we will stretch the concept of an algorithm to also allow for inclusion of not precisely defined subroutines, such as `pick next search node to investigate'. As will become apparent, B\&B contains many such subroutines that need to be described in a precise way in order to make it an implementable algorithm. 

B\&B systematically explores the feasible region by partitioning it into sub-MILPs and obtaining bounds via their LP relaxations, which serve as a mechanism to rule out suboptimal regions. This partitioning scheme can be represented by a search tree, where each node represents a sub-MILP. We use the notation $z^{LP_i}$ to denote the optimal value to the LP relaxation of node $i$, and we denote $\vx^{LP_i}$ its optimal solution. At any point during the search, the value of an integer feasible solution provides an upper bound on $z^*$. The best, i.e., smallest, known upper bound value is called \emph{primal bound}, and we denote it by $\bar{z}$. Similarly, the bounds $z^{LP_i}$ of unprocessed nodes can be used to obtain a lower bound. In particular, $\underline{z}:=\min_{i: i\text{ unprocessed}} \{z^{LP_i}\}$ provides a lower bound on $z^*$ and is called \emph{dual bound}. The B\&B algorithm ends when all nodes have been processed, when $\bar{z}=\underline{z}$, or when another termination condition (e.g., timeout) applies. 

These are the 
basic principles of the  
B\&B algorithm as implemented in commercial solvers, such as CPLEX \cite{CPLEX}, Gurobi \cite{gurobi} or Xpress \cite{xpress}, as well as in academic solvers like SCIP \cite{SCIP7}. In practice, the execution of the B\&B algorithm revolves around some key solver components that handle the different aspects of the solving process. The most important components are preprocessing, the branching rule, the cut management and the primal heuristics.

\paragraph{Preprocessing.} Most solvers implement a number of procedures that try to reduce the size of the problem and its difficulty, for example by identifying substructures or ways to strengthen the LP relaxation.

\paragraph{Branching and node selection.} The procedure of dividing the feasible region is called \emph{branching}. There is a choice to be made with respect to which disjunction is used for branching. The standard is to use single variable disjunctions of the type
$$x_j \leq \lfloor x^{LP}_j \rfloor \lor x_j \geq \lceil x^{LP}_j \rceil\,, $$
for some $j\in\mathcal{I}$ such that $x^{LP}_j \notin \mathbb{Z}$. Still, the solution to the LP relaxation is likely to violate more than one integrality constraint, which means there is more than one candidate variable for branching. The so-called \emph{branching rule} is the strategy that the solver uses to select a variable for branching. Computational studies have shown that this choice has a critical impact on solver performance \cite{Achterberg2013}. Another important decision, though less critical in terms of solver performance, is \emph{node selection}, where the question is which unprocessed B\&B node to consider next.

\paragraph{Cutting planes.} The LP relaxation can be strengthened by adding valid linear inequalities. These are inequalities that cut off parts of the relaxation but do not exclude any integer feasible solutions. This can, in principle, be done in any node of the B\&B tree representing a non-empty feasible region, yet it is standard practice to use cutting planes more heavily, or even exclusively, at the root node. A B\&B routine where cutting planes are added in other nodes than the root node is typically referred to as branch-and-cut \cite{PadRin1991}. There is a vast amount of knowledge on cutting planes (or cuts, for short) for MILP (see, e.g., \cite{NemhauserWolsey, Wolsey, Conforti2014}) and most solvers implement a plethora of efficient separation\footnote{The term \emph{separation} refers to the fact that a cutting plane has an effect when it separates a fractional solution of an LP relaxation from the convex hull of (mixed) integer feasible solutions.} algorithms. The wide availability of cuts can create an issue: while the addition of cuts strengthens the local relaxation, a large amount of cuts can slow down LP solving and lead to numerical instability. For this reason, a judicious \emph{cut management} strategy, which includes separation, selection and removal, is of utmost importance. 

\paragraph{Primal heuristics.}
We use the term primal heuristics to refer to routines that try to find feasible solutions in a short amount of time without a success guarantee of doing so. Relying solely on integer feasible LP relaxations to find solutions is most often inefficient. Primal heuristics can provide good solutions early on and help bring down the primal bound $\bar{z}$ more effectively. As with cutting planes, primal heuristics can be used in any node of the tree if desired.

\subsection{Evaluation metrics for MILP}
\label{subsec:metrics}
In this section, we define a number of metrics that quantify the progress of a branch-and-bound run. We extend the notation just presented in Section~\ref{subsec:intro_to_milp} with a variable $t\geq 0$ that represents solving time. We can then define the primal bound $\bar{z}(t):[0,T_{\max}]\mapsto \mathbb{R}$ and the dual bound $\underline{z}(t):[0,T_{\max}]\mapsto \mathbb{R}$ as functions of time. We define $\bar{z}(t)$ to be infinity if no integer feasible solution has yet been found at time $t$. We also make use of a small tolerance value $\epsilon$ that is introduced for numerical stability, typically $\epsilon = 10^{-6}$. Notice that many of these metrics make use of the optimal solution value $z^*$ and must therefore be calculated a posteriori, once the instance is solved. 

\subsubsection*{Optimality gap}
We follow the definition of SCIP\footnote{\url{https://www.scipopt.org/doc/html/group__PublicSolvingStatsMethods.php}} 
of (normalised) optimality gap, namely
\begin{equation}
\label{eq:gap}
    g(t) :=
    \begin{cases*}
      \infty & if no solution has been 
      found yet or $\bar{z}(t)\cdot \underline{z}(t)<0$,\\
      \frac{|\bar{z}(t)-\underline{z}(t)|}{\min \{|\bar{z}(t)|, |\underline{z}(t)|\}} & otherwise.
\end{cases*}
\end{equation}

Alternatively, one can track the \emph{integrality gap} defined as $g'(t)=|z^*-\underline{z}(t)|$.

\subsubsection*{Primal gap and integral}
For a given feasible solution $\vx$, we define the \emph{primal gap} $\gamma(\vx)$ as
\begin{equation}
\label{eq:primalgap}
   \gamma(\vx) :=
    \begin{cases*}
      1 & if $z^*\cdot \vc^T\vx < 0$, \\
      \frac{|z^*-\vc^T\vx|}{\max\{ |z^*|,|\vc^T\vx|, \epsilon \}}        & otherwise.
    \end{cases*}
  \end{equation}

  We can define a \emph{primal gap function} that maps the solving time to the primal gap of the best solution found up until that point. In particular, denoting $\vx(t)$ the best solution found at time $t$, we define
  \begin{equation}
\label{eq:primalgapfun}
   p(t) :=
    \begin{cases*}
      1 & if no solution has been found at time $t$,\\
      \gamma(\vx(t)) & otherwise.
    \end{cases*}
  \end{equation}

  The \emph{primal integral} \cite{Berthold2013} of a process with time limit $t_{max}>0$ is defined as
  \begin{equation}
  \label{eq:primalintegral}
      P(t_{max}) := \int_0^{t_{max}} p(t) dt.
  \end{equation}

  \subsubsection*{Primal-dual integral}
One can extend the concept of primal integral to also account for improvements in the dual bound. For this purpose, we use the optimality gap instead of the primal gap, and integrate the function
      \begin{equation}
    \label{eq:primaldualfun}
   pd(t) :=
    \begin{cases*}
      1 & if no solution has been 
      found yet or $\bar{z}(t)\cdot \underline{z}(t)<0$,\\
    g(t)  & otherwise,
    \end{cases*}
  \end{equation}
to obtain the primal-dual integral, defined as
    \begin{equation}
  \label{eq:primaldualintegral}
      PD(t_{max}) := \int_0^{t_{max}} pd(t) dt.
  \end{equation}

\section{A brief introduction to machine learning}
\label{sec:into_to_ml}
This section provides a brief introduction to the key 
concepts of machine learning that are necessary in order to follow this survey. 
For a more detailed introduction to machine learning we refer to \citet{Mitchell17:book}.

We are interested in the problem of constructing a mapping from some input data to a desired output space. Let $\mathcal{X}$ be the input data space and let $\mathcal{Y}$ be the output space. It is common to refer to $\mathcal{X}$ as the \emph{feature} space, as it basically represents a set of descriptors of the data samples. The output can be, for example, a prediction based on the input data. The mapping will be constructed by choosing among a family of parameterized functions $f(\cdot, \vtheta): \mathcal{X} \to \mathcal{Y}$ with parameters $\vtheta \in \Theta$. 
In short, the objective is to optimize the behaviour of the mapping $f(\cdot, \vtheta)$ by carefully tuning the parameters, based on sampled observations and a progress metric of choice.

To formalize this, we follow standard practices in machine learning and distinguish the following two settings. 

\paragraph{Supervised learning.} 
The learner has access to a finite collection of pairs \linebreak[4]$\{(X_i, Y_i)\}_{i=1}^N$, where $X_i\in\mathcal{X}$ and $Y_i\in\mathcal{Y}$. Here, $Y_i$ is the \emph{desired output} to input $X_i$. It is common to refer to $Y_i$ as the \emph{ground truth} or \emph{label}. The goal is to minimize the \emph{loss function} $l:\mathcal{Y}\times \mathcal{Y} \mapsto \mathbb{R}$, a metric that represents the dissimilarity between the prediction and the desired output. An example is the mean square error loss (MSE) $\frac{1}{N}\sum_{i=1}^N \Big(f(X_i, \vtheta) - Y_i\Big)^2$. Since the true sample distribution is unknown, the loss is minimized with respect to the observed samples
$$\underset{\theta \in \Theta}{\text{minimize}} \sum_{i=1}^N l(Y_i, f(X_i,\theta))\,.$$

\begin{figure}[tb]
    \centering
    \includegraphics[width=0.7\textwidth]{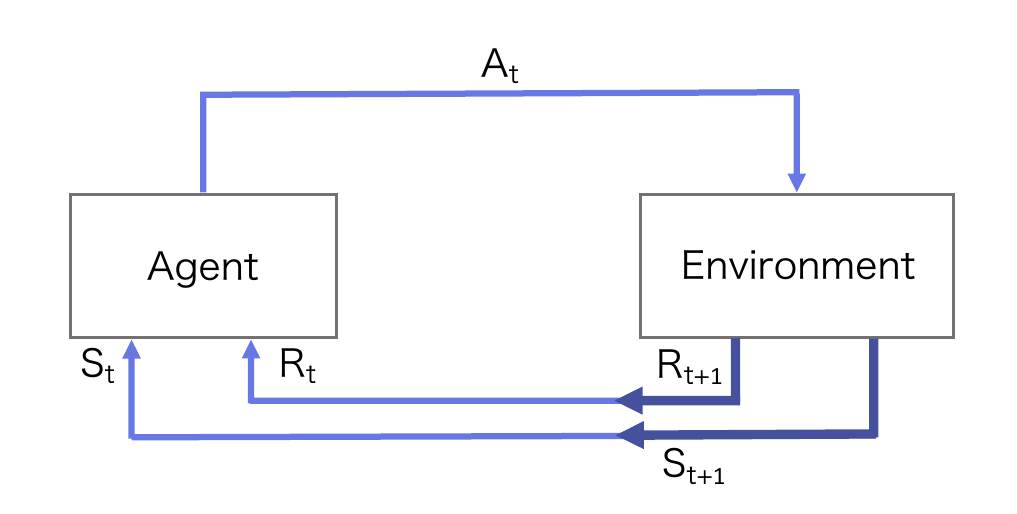}
    \caption{Markov Decision Process}
    \label{fig:mdp}
\end{figure}

\paragraph{Reinforcement learning.}
Reinforcement learning (RL) is defined in the context of sequential decision making, where actions have long-term consequences and the optimal action is either unknown or too expensive to compute. 
The learning agent has no access to the ground truth. This methodology is formalized under the framework of Markov Decision Processes (MDPs), see Figure~\ref{fig:mdp} for a diagram summarizing MDPs. In an MDP, an \emph{agent} interacts with an \emph{environment}. The environment has an associated \emph{state} representing its internal configuration. We denote the state space $\mathcal{S}$. The agent acts on the environment by choosing an action from the action space $\mathcal{A}$ using a \emph{policy} $\pi:\mathcal{S}\mapsto \mathcal{A}$. The agent's action changes the state of the environment, which corresponds to the transition to a new state. The transitioning mechanism is unknown to the agent.

Apart from the environment's state, the agent can observe a \termdefn{reward function}. We consider the episodic case, where this interaction between agent and environment happens sequentially in discrete time-steps until a \emph{terminal state} is reached and the interaction ends. A realization of such agent-environment interaction is called trajectory
$$\tau := (S_0, A_0, R_1, \dots, S_{T-1}, A_{T-1}, R_T, S_T)\,,$$
where $T$ denotes the episode length. The goal of the agent is to find a policy that maximizes the expected cumulative reward, known as the \termdefn{value function} and formally defined as
\begin{equation}
\label{eq:rl}
    V_\pi := \mathbb{E}_{\tau} \Biggl[ \sum_{t=1}^T \gamma^{t-1} R_t \Biggr]\,,
\end{equation}
where $\gamma \in [0,1]$ is the \emph{discount factor}. This parameter controls the greediness of the policy: for $\gamma$ close to zero, the agent will prioritize obtaining immediate rewards, whereas for $\gamma$ close to one the agent is encouraged to follow a strategy that pays off in the long term. There is a variety of RL methodologies (see, e.g., \cite{Sewak2019}) that provide ways to train a parametrized function $\pi(s, \theta)$ with the goal of maximizing $V_\pi$. Notice that the trajectory distribution depends both on the agent's policy and the unknown transition mechanism of the environment. For this reason, RL algorithms use trajectory sampling as a way to estimate this expectation. Note that $-V_\pi$ can be seen as a loss function.

The MDP formulation can also be used in a supervised learning setting. In this case, the reward signal is substituted by an expert that tells the agent what is the optimal action.  The expert is typically expensive to query, and for this reason we want the agent to learn a policy that imitates the expert, but at a lower cost. This methodology is known as \emph{imitation learning}.


\subsection{Mapping features to predictions}
\label{subsec:nn}
A fundamental step in the design of a ML scheme is the choice of function space parameterized by $\vtheta$. This is commonly known as the \termdefn{architecture}, and it is independent of the choice of learning methodology. Modern machine learning favors neural networks as the computational representation and mechanism of $f$. A simple form of neural network is the feed-forward neural network that consists of a series of linear transformations, each followed by a non-linear transformation called an \termdefn{activation function}.


\begin{definition}[Feed-forward Neural Network]
    Let $X\in \mathbb{R}^d$ be the input data and let $d^L$ be the desired output dimension. An $L$-layer \emph{feed-forward neural network} is a function $f:\mathbb{R}^d \mapsto \mathbb{R}^{d^L}$ that defines a mapping from $X$ to $a^L \in \mathbb{R}^{d^L}$ through the recursive relation
    $$z^{l} = W^{l} a^{l-1} + b^{l}$$
    $$a^{l} = \sigma (z^l)$$
    for $l=1,...,L$, where 
    \begin{itemize}
        \item $a^0=X$, $d^0=d$,
        \item $W^l \in \mathbb{R}^{d^l \times d^{l-1}}$ and $b^l\in \mathbb{R}^{d^l}$ are learnable parameters, and
        \item $\sigma$ is the activation function.
\end{itemize}
\end{definition}

Some common activation functions are the Rectified Linear Unit (ReLU) or the sigmoid function, defined, respectively, as
$$\sigma_{\text{ReLu}}(x) = \max (0, x)\,,$$
$$\sigma_{s}(x) = \frac{1}{1+e^{-x}}\, .$$
Notice that these functions are applied componentwise.

Another prominent architecture is the graph neural network (GNN) that will be of particular interest to the discussion in Section~\ref{sec:representations}. Consider an undirected graph $G=(V,E)$ with vertex set $V$ and edge set $E$. For a node $v\in V$, we denote $N(v)\subseteq V$ the set of neighbors of $v$.

\begin{definition}[Graph embedding]
    A ($d$-dimensional) \emph{graph embedding} $\xi$ is a function that takes in a graph $G=(V,E)$ and a node $v\in V$ and returns an element $\xi (G,v) \in \mathbb{R}^d$.
\end{definition}

\begin{definition}[Graph Neural Network]
    \label{def:gnn}
    A \emph{Graph Neural Network} is a function that takes as input a graph $G=(V,E)$ and an initial embedding $\xi^0$ and defines a recursive embedding $\xi^{t}$ over the vertices of $G$. This function is characterized by
    \begin{itemize}
        \item A combination function \texttt{comb}: $\mathbb{R}^{2d}\to \mathbb{R}^d$, and
        \item An update rule $\xi^{t+1}(G,v) = \texttt{comb} \left(\xi^{t}(G,v), \sum_{u\in N(v)} \xi^{t}(G,u) \right)$.
\end{itemize}
\end{definition}

\begin{figure}
    \centering
    \includegraphics[width=0.7\textwidth]{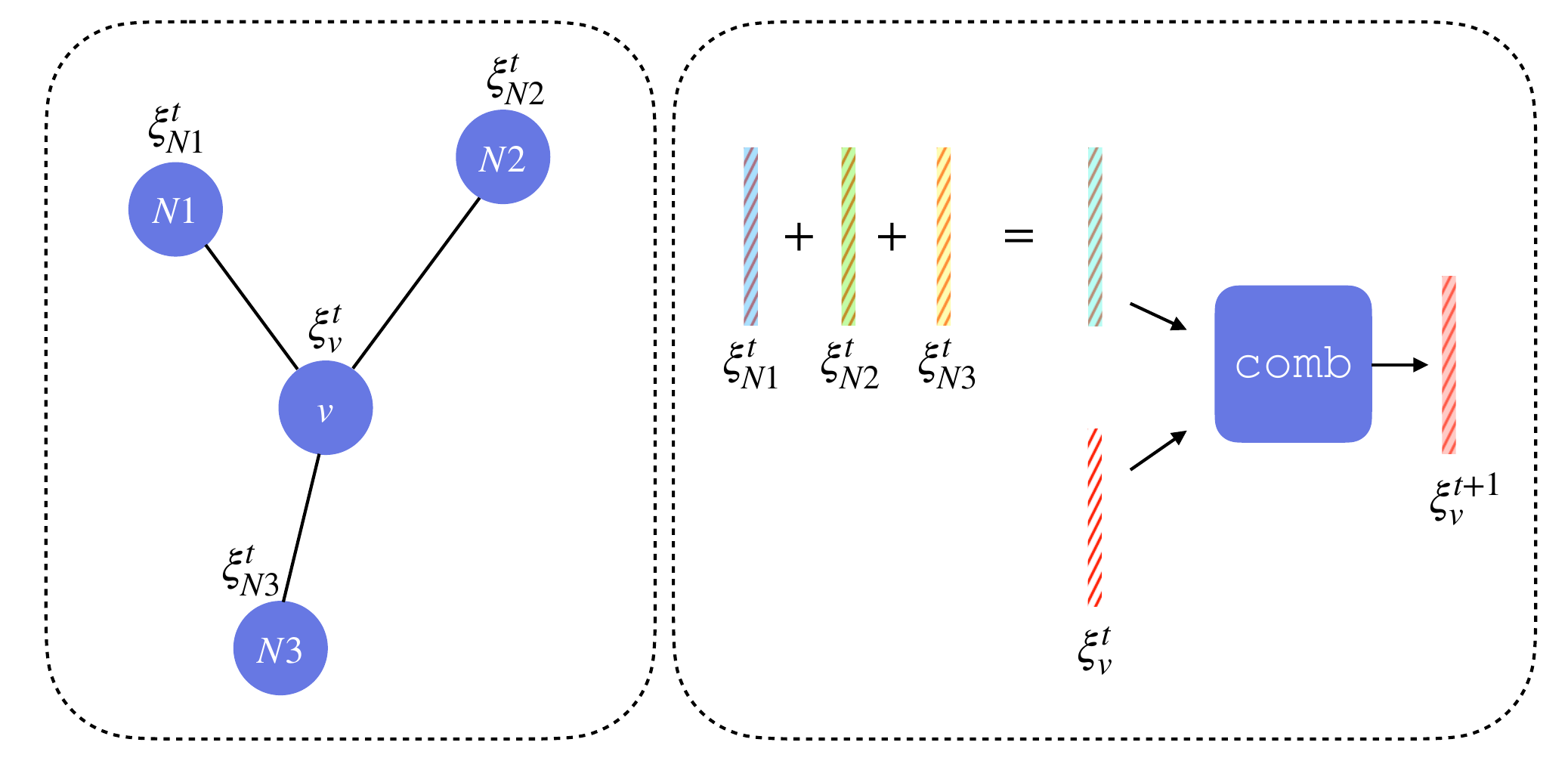}
    \caption{Embedding computation with a Graph Neural Network. Here, we use the abbreviation $\xi^t_v$ for $\xi^t(G,v)$. To update the embedding $\xi^t_v$ of node $v$ at time $t$, the embeddings of neighboring nodes are added and then combined with the current embedding using the \texttt{comb} function. The result is $\xi^{t+1}_v$.}
    \label{fig:gnn}
\end{figure}

Figure \ref{fig:gnn} depicts one iteration of this recursive embedding mechanism. It is common to refer to a single iteration as a \emph{message-passing operation}. A possible combination function is a feed-forward neural network over the two concatenated vectors. We can extend graph embeddings to edges. This is, for $e\in E$, we can define $\xi (G,e)\in\mathbb{R}^d$. In the presence of edge embeddings, we can redefine a GNN update rule as
\begin{equation}
\label{eq:gnns}
    \xi^{t+1}(G,v) = \texttt{comb} \left(\xi^{t}(G,v), \sum_{u\in N(v)} \texttt{aggr} \left( \xi^{t}(G,u), \xi(G,\{u,v\}) \right) \right)
\end{equation}
with an aggregation function \texttt{aggr}: $\mathbb{R}^{2d}\to \mathbb{R}^d$ that handles the combination of node and edge embeddings. An example of aggregation function is a simple component-wise sum.

\subsection{Elements of the learning process }

The learning process is itself one of (inexact) optimization, termed \termdefn{training}.

\paragraph{Optimization algorithms.} The goal of the training process is to optimize the behaviour of the mapping $f(\cdot, \vtheta)$ over the parameter space $\Theta$ with respect to the loss function. This is done with algorithms for continuous non-convex optimization, which take gradient descent as their base.  Some examples of commonly-used algorithms are AdaGrad \citep{Duchi2011} or Adam \citep{Kingma2014}. For reinforcement learning, the situation is more complex, in the sense that several learning paradigms are available for approximating the value function and learning the optimal policy.\footnote{A general discussion of those paradigms is outside the scope of the paper and special cases will be discussed while surveying specific results from the literature, see Section \ref{sec:tasks}.} 

\paragraph{Hyper-parameters.} There is a number of additional parameters that influence the learning process that are not part of the parameters $\vtheta$ to be learned. These meta-level parameters are called \termdefn{hyper-parameters}. Examples of hyper-parameters are the number of iterations of the optimization algorithm, the learning rate, the size of the learned function, etc. Tools exist to automatically tune hyper-parameters given a validation dataset (see below).

\paragraph{Train, validation and test datasets.}
The data used to tune the parameters of mapping $f$ is called the \emph{training data}. The data used to estimate the performance of $f$ in an unbiased way, and to tune the hyper-parameters, is called the \emph{validation data}. Finally, the data used to evaluate the final trained and optimised mapping is called the \emph{test set}. These sets must be disjoint.

\paragraph{Data collection.} The learner makes use of the available training data to improve its behaviour. This data can be gathered all at once during an initial data collection phase, before training starts. Alternatively, several data collection and training steps can be run in an alternating fashion.

\paragraph{Overfitting.} We say a learned model suffers from \emph{overfitting} when the performance on the training set is distinctively better than the performance on the test set. In other words, the model is unable to generalise to data unseen during the training phase. This phenomenon occurs when the chosen model is too large (in terms of number of parameters) for the task, or the amount of relevant data is too limited for the size of the model being learned. See \citet{Goodfellow2016} for a discussion of methods to avoid overfitting.

\paragraph{Online versus offline learning.} Whenever we can distinguish separate training and execution phases we speak of \termdefn{offline learning}.  That is, in offline learning, we perform the data collection and learning as a separate, preliminary, once and for all process.  After that, the learned function is fixed and used without further tuning. This training phase can be computationally costly but this cost is not reflected in the execution time (often called inference time) of the algorithm, which is typically the metric of interest to the end user. In fact, this can be seen as doing the heavy work upfront while alleviating the effort at the moment of execution.

In contrast, in an \termdefn{online learning} setting, data becomes available during execution. Learning must therefore happen dynamically as each data sample becomes available, and in parallel to the main process. This has the great advantage of generating a function that is highly adaptive. On the other hand, it entails adding the cost of learning to the cost of execution and further there is usually less available data.

\section{Learning tasks}
\label{sec:tasks}

So far, we have introduced the elements of the B\&B algorithm for MILP, and introduced general terminology in machine learning. This section presents different ways to formulate a learning task within the B\&B algorithm. We split this discussion by considering each major solver component separately. 

Before getting into the details, two considerations are required. First, throughout the section we will use the symbol $X$ to abstractly denote a representation of an MILP instance. This representation may include data coming from the problem description (such as in Equation~\ref{eq:MILP}) and from the B\&B process. Section~\ref{sec:representations} discusses, in more detail, approaches to build such a representation. Second, it is important to highlight that the learning tasks that we describe next might be associated with different levels of required generalization. More precisely, 
the ML models are often trained on data belonging to MILP instances in the same class, for example, set covering, knapsack, etc. In other words, the characteristics of the constraint matrix in \eqref{eq:MILP} are leveraged to obtain accurate predictions, so the generalization power of the resulting models is generally restricted to that specific MILP class used for training. However, we will also review cases in which learning  happens instance by instance or extends outside of a known MILP distribution. We will highlight these differences throughout the entire review of the learning tasks.

\subsection{Primal heuristics}
\label{subsec:primal}

\begin{figure}[tb]
     \centering
     \includegraphics[scale=0.72]{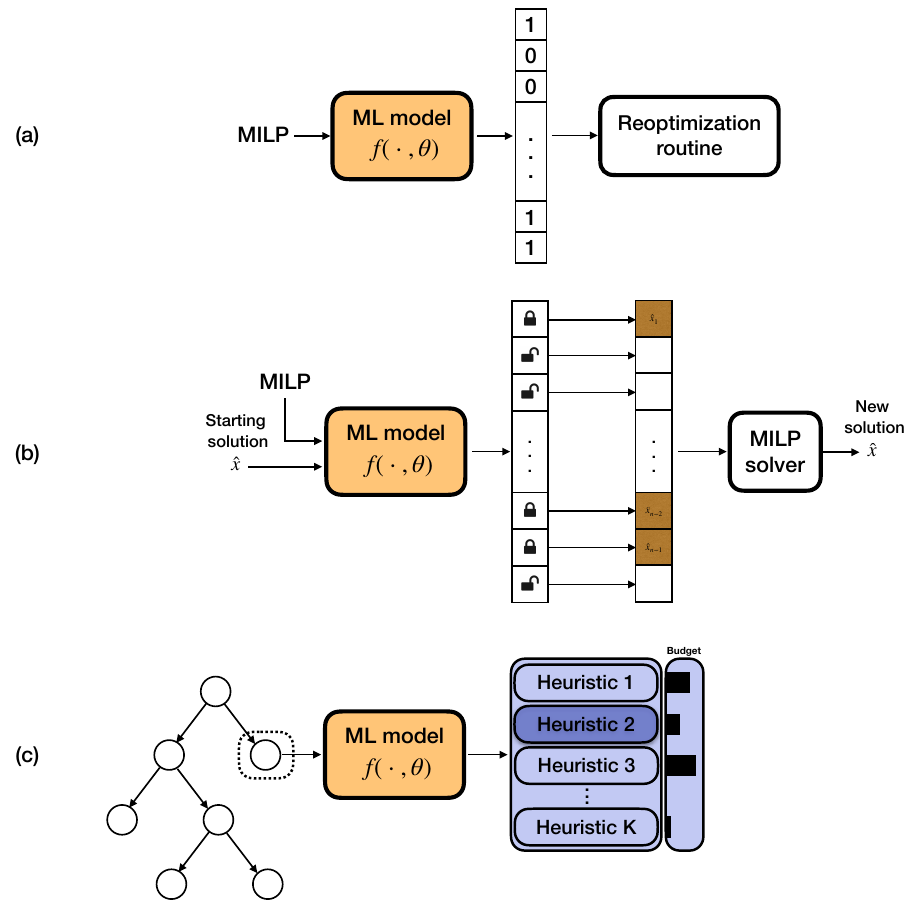}
     \caption{Three learning problems related to primal heuristics: (a) predict a reference solution and search in its neighborhood, (b) neighborhood selection -- which and/or how many variables to unfix and re-optimize, (c) heuristic scheduling -- which heuristics to run and/or for how long. }
     \label{fig:primal_models}
 \end{figure}

Primal heuristics play a crucial role in quickly finding feasible solutions and consequently improving the primal bound $\Bar{z}$ 
in the early stages of the solve. In the context of primal heuristics, ML techniques can be of interest to leverage common structures in the instances. A number of methodologies have been proposed for this purpose. Conceptually, they can be split into three main approaches: (i) guiding a heuristic search with a starting predicted solution, (ii) solution improvement via a learned 
neighborhood selection criterion, and (iii) learning a schedule to pre-existing heuristic routines. These methodologies are summarized in Figure~\ref{fig:primal_models} and will be discussed in more detail in this section.

Typically, we divide primal heuristics into two groups. We speak of \emph{improvement heuristics} when referring to routines that need a feasible solution to begin with, and whose goal is to improve upon it. Otherwise, we speak of \emph{starting heuristics}, which can be run at any moment during the search process. A prime example of improvement heuristics is Large 
Neighborhood Search (LNS) \cite{PisRop2019}. The idea is to optimize an auxiliary MILP of smaller size, constructed by reducing the feasible region of the original MILP. This is typically done by fixing the value of some of the variables and optimizing the rest. Another strategy is to search over the neighborhood of a solution $\hat{\vx}$ by imposing the constraint
\begin{equation}
\label{eq:local_branching}
    \sum_{j:\hat{\vx}_j=0}x_j + \sum_{j:\hat{\vx}_j=1}(1-x_j) \leq \eta\, ,
\end{equation}
with $\eta$ the parameter that controls the neighborhood size. This is known as \emph{local branching} \cite{Fischetti2003}.

\subsubsection{Solution prediction to guide the search} 
 Some authors explore the idea of using predictions of the optimal solution. The goal is to produce a (partial) assignment of the binary variables\footnote{Note that the focus on binary variables is justified because, in most well-established instance collections, binary variables account for the great majority of integer variables. As an example, more than 68\% of the instances in the MIPLIB \cite{MIPLIB} are purely binary, and in the remaining ones, more than 90\% of the integer variables are in fact binary.  See, e.g., \citet{Nair2020} for extensions to integer variables.} in a binary or mixed binary MILP, that can then be used to guide the search. This idea has been implemented in different ways, both in terms of how to obtain the predictions and how to use them. Let us start by discussing the latter.

 \citet{Ding2020} use a local branching constraint (see Eq.~\ref{eq:local_branching}) with respect to predicted values of a subset $J\subseteq \mathcal{B}$ of the binary variables. This restricts the search to a neighborhood of the predicted partial solution. In contrast, both \citet{Nair2020} and \citet{Khalil2022mipgnn} propose to fix the variables in $J$ to their predicted value and hand over this partial assignment to an 
 MILP solver that optimizes over the remaining variables.  In MILP terminology, this corresponds to a warm start.  \citet{Khalil2022mipgnn} further use the predictions to guide node selection (see Section \ref{subsec:nodesel}).
 
The question still remains as to how to obtain these predictions. This problem is naturally formulated as a supervised learning task, where the desired output is the optimal solution. However, optimal solutions to be used as labels in the training phase are costly to obtain and do not capture information about the region where they lie. The aforementioned work of \citet{Ding2020, Nair2020, Khalil2022mipgnn} make use of a set of feasible (not necessarily optimal) solutions to learn predictions. This is, the learning process starts with a data collection phase that, for each instance $X$, gathers a set of feasible solutions $D(X) = \{\hat{\vx}^{(1)},\hat{\vx}^{(2)},..., \hat{\vx}^{(K)}\}$. The goal is to learn a function $p_\theta (x_j=1|X)$, with parameters $\theta$, that can be interpreted as the probability that variable $x_j$ takes value 1 given problem $X$. The parameters are tuned to make the behavior of these functions resemble as closely as possible that of a target probability distribution $p_T (x_j=1|X)$. Table \ref{tab:primal_target} summarizes the proposed targets.

 \begin{table}[tb]
\centering
\caption{Three methods to define the target probability distribution $p_T$, given a collection $D(X) = \{\hat{\vx}^{(1)},\hat{\vx}^{(2)},..., \hat{\vx}^{(K)}\}$ of feasible solutions to problem $X$. 
}
\label{tab:primal_target}
\begin{tabular}{lc}
Authors  & Target   \\ \hline
\multirow{2}{*}{\citet{Ding2020}} & $p_T (x_j=1|X) = \mathds{1}_{\hat{x}_j^{(1)}=1}$ for every $j$ s.t. \\
 &  $\hat{x}_j^{(k)}=\hat{x}_j^{(1)}$ for all $k\in \{1,...,K\}$ \\
 & \\
\multirow{2}{*}{\citet{Nair2020}} & $p_T (\hat{\vx}^{(k)}|X) = e^{-\vc^T \hat{\vx}^{(k)}} / \sum_{i=1}^K e^{- \vc^T \hat{\vx}^{(i)}}$ \\
 &  and $p_\theta (\hat{\vx}|X) := \prod_{j=1}^n p_\theta (x_j=\hat{x}_j|X)$  \\
 & \\
 \citet{Khalil2022mipgnn} & $p_T (x_j=1|X) = \frac{1}{K}\sum_{k=1}^K \hat{x}^{(k)}_j$ \\
\end{tabular}
\end{table}

\subsubsection{Solution improvement via neighborhood selection}
Alternatively to solution prediction, one may be interested in learning a \textit{destroy heuristic} criterion. That is, given an initial feasible solution $\hat{\vx}$, we select a subset of the integer variables to be re-optimized, leaving the remaining integer variables fixed. This process can be run iteratively. The goal is to identify substructures of the problem that can be used to decompose it into smaller, more manageable sub-problems.

Following this line, \citet{Song2020} learn to partition the set of integer variables $\mathcal{I}$ into $K$ disjoint subsets $\{S_i\}_{i=1}^K$ such that $\mathcal{I} = S_1\cup ...\cup S_K$. They iteratively unfix and reoptimize the variables in each subset, fixing the rest to the value in the best known solution. The variable subsets are of fixed size, which means that the variables simply need to be classified into the subset they belong to. Alternatively, \citet{Wu2021}, \citet{Sonnerat2021} and 
\citet{Huang2023} propose a method to select a single flexible-sized subset of variables to unfix. Their approach resembles the value-prediction methodology discussed in the previous section. However, instead of predicting the value that a variable takes in an optimal solution, they aim at predicting whether the variable is already assigned to its optimal value in the current best solution $\hat{\vx}$. The number of unfixings is learnt implicitly through the conditionally-independent probabilities $p_\theta (\hat{x}_j = x^*_j|X)$ from which they sample variable unfixings. While \citet{Wu2021} use a reinforcement learning algorithm to train their policy, \citet{Sonnerat2021} explicitly calculate the best solution at most $\eta$ unfixings away from $\hat{\vx}$ and use it as a target unfixing policy. \citet{Huang2023} follow a similar approach, while also providing the learner with negative examples. This means variable unfixings that do not lead to a sufficiently large improvement, as compared to the best. This provides additional information for the learner to distinguish good and bad unfixings. 

The methods of \citet{Song2020, Wu2021, Sonnerat2021} and \citet{Huang2023} were compared against continuously running the 
MILP solvers they use as a subroutine. Given the same amount of time, the ML-assisted methods are able to find better solutions. Interestingly, using a random selection of the variables to unfix often results in a performance improvement compared to running the solver continuously. This highlights the fact that MILP solvers are typically tuned to optimize different performance metrics, such as the optimality gap, which means that a significant part of the computational effort is spent in, e.g., obtaining good dual bounds. Nonetheless, the proposed methods provide a relevant methodology to make use of the MILP solver in a black-box manner when the goal is to obtain high-quality but not necessarily optimal solutions in a short amount of time.

Related to this line of work, \citet{Liu2022} point out that the optimal value of the neighborhood size parameter $\eta$ is strongly dependent on the class of instances being solved. They work in the context of local branching and propose to learn two policies. A first policy $f_{init}(X)$, obtained in a supervised way, determines an appropriate neighborhood size $\eta_0$ for a first iteration. Another policy $f_a(X)$, obtained by RL, decides how to \emph{adapt} the neighborhood size in successive iterations given the previous solver statistics. Compared to the simple neighborhood size selection proposed in \citet{Fischetti2003}, these two policies close the primal gap (see Eq.~\ref{eq:primalgap}) faster, especially when combined together. Their computational experiments demonstrate that these selection rules have great generalization ability across different instance sizes and types, especially for $f_a(X)$, also across different classes of instances, proving the wide applicability of their approach.

\subsubsection{Learning to schedule heuristics}
Whether they use ML or not, many successful primal heuristics for MILP have been proposed. Experimentally, we observe that no single heuristic dominates the others on all problems \cite{Hendel2019}. Their performance is highly dependent on the problem, and even on the solving stage. A new decision problem arises: given a collection of primal heuristics, which one should be run? In a broader sense, heuristic scheduling also tries to answer questions such as how frequently should heuristics be run, or under what computational budget. In this section, we review the line of research on using learning methods to answer such questions.

\citet{Hendel2022} proposes Adaptive Large Neighborhood Search (ALNS), a heuristic that, whenever called, makes a choice among eight Large Neighborhood Search methods. This choice is framed as a Multi Armed Bandit problem \cite{Bubeck2012}, an online learning methodology that learns a selection policy per instance. The learned policy is based on the observed (a posteriori) performance of the chosen heuristics. This simple formulation encapsulates the classical exploitation versus exploration dilemma: we must balance running heuristics that have performed well in the past with running heuristics whose performance is unknown because they have not been selected a sufficient number of times. The author uses the $\alpha$-UCB algorithm \cite{Bubeck2012} with a reward function that combines several aspects of a heuristic's performance: whether an incumbent was found, the objective improvement and the computational time the heuristic needs. The experiments on MIPLIB 2017 \cite{Gleixner2021} show a considerable improvement in the primal integral (see Eq.~\ref{eq:primalintegral}), as well as a speed-up. The Multi Armed Bandit formulation was also used by \citet{Hendel2019} and \citet{Chmiela2023}, who extend it to new types of heuristics.  

\citet{Chmiela2021} take a different perspective: instead of a selection policy per instance, they create a prioritization order that applies to a given class of instances. In fact, these authors' method constructs a schedule that assigns both a priority and a computational budget to each heuristic. This schedule is crafted offline, after a data collection phase where different heuristics are evaluated. The objective is to minimize the computational budget while finding feasible solutions for a large fraction of the B\&B nodes. 
 The proposed method produces schedules for a number of different heuristics, which result in improvements in the primal integral compared to default SCIP and even a version of the solver that has been manually tuned for the considered class of instances.

Other than deciding which heuristics to run and for how long, there is the question of \emph{when} to run them. The work discussed so far uses conventional rules to decide at which nodes heuristics will be run. The methodology proposed by \citet{Khalil2017} is to instead build a mapping between B\&B nodes and a yes/no decision for running each heuristic. 
 This decision problem is challenging. Even with perfect knowledge of when a heuristic will be successful, a run-when-successful rule does not necessarily minimize solving time. \citet{Khalil2017} analyze the competitive ratio of such a rule compared to an optimal offline decision policy. 
 They then imitate the imperfect run-when-successful rule by learning to predict when a heuristic will succeed. 
 Their computational study shows the effectiveness of this method: heuristics are called less often, but with a higher success rate, resulting in a better primal integral. This effect is even stronger when the policy is trained on data coming from the same instance class.

\subsubsection{Future outlook}
The success of most of the approaches discussed in Section \ref{subsec:primal} is rooted in the ability of learning within the distribution of specific classes of MILP instances. Generalizing outside of a specific 
class, has proven so far difficult. This is certainly the main challenge for the integration of ML-augmented primal heuristics within MILP solvers. Studies like the one in \citet{Liu2022} indicate that some significant generalization is achievable, especially in the context of algorithms that sequentially adapt while exploring the solution space as in classical RL schemes. This is a promising direction that, however, conflicts with the MILP solvers general need of executing (extremely) fast primal heuristics.


\subsection{Branching}
\label{subsec:branching}

Branching is one of the core mechanisms on which the B\&B algorithm operates. The branching rule, i.e., the criterion used to select a variable for branching, has been identified as having a critical impact on performance \cite{Achterberg2013}. Computational studies have served to identify a number of metrics that are good indicators of how a variable will perform. Notably, state-of-the-art branching rules look at the change in objective value in the resulting children nodes. 
More specifically, let $z^{LP_i}$ be the objective value of the current node $i$ and let $z^{LP_{i+}}$ and  $z^{LP_{i-}}$ be the objective values of the two nodes resulting from branching on candidate variable $x_k$. This variable will be scored using a combination\footnote{For example, the variable's score can be computed as $\max\{\Delta^+, 10^{-8}\}\cdot \max\{\Delta^-, 10^{-8}\}$. See, e.g. \citet{Achterberg2007} for a discussion.} of $\Delta^+:=z^{LP_{i+}}-z^{LP_{i}}$ and $\Delta^-:=z^{LP_{i-}}-z^{LP_{i}}$. These values can be explicitly calculated for each candidate at the moment of branching (thus solving two LPs per candidate). Such strategy, known as \emph{strong branching}, was introduced in the context of the travelling salesman problem \cite{Applegate1995} and later standardised by CPLEX. 

Strong branching is known to produce small B\&B trees, but at high computational cost per branching. Alternatively, one can attempt to estimate the objective change based on past values, once they become naturally available through node processing. In particular, solvers typically store the values of $\Delta^+$ and $\Delta^-$ normalized by the variable's fractionality, and keep track of the per-variable average, known as \emph{pseudocosts}. The so-called \emph{reliability branching} rule \cite{Achterberg2005} performs strong branching at the top of the tree as an initialization phase, and then switches to using pseudocosts as soon as a variable has been branched on enough times. The initialization phase is not only important to build a branching history, but also because branching decisions have the most impact at the top of the tree. 
This is because judicious branching decisions here can lead to much smaller tree sizes, due to the earlier finding of feasible solutions and the stronger pruning of nodes. 

In this section, we will discuss different approaches to learning to branch. Their common goal is to learn a function that maps a description of the candidate variables to scores. Ultimately, this is done with the objective of minimizing solving time, which typically entails a favorable balance among different sub-targets. For example, (i) being computationally cheap, (ii) generating small trees as a result of their scoring, and (iii) adapting to the different situations that may arise. These objectives are often at odds with each other.

\subsubsection{A first approach to learning from strong branching}
There is a well-established body of research on 
fast approximations of the strong branching rule. This idea was first explored by \citet{Alvarez2017} who propose to learn a prediction of the strong branching score of each variable. The predictor is learned offline (see Section~\ref{sec:into_to_ml}) and tested on both heterogeneous and homogeneous instance collections. The experiments of \citet{Alvarez2017} show that the method achieves the desired objective of imitating strong branching decisions without the large computational overhead. Indeed, when compared with strong branching on a fixed number of nodes, the closed gap is only moderately worse, indicating that the branching decisions are of high quality. Simultaneously, for a fixed time, the method explores a much larger number of nodes and closes a greater proportion of the gap, therefore achieving an overall better trade-off. However, reliability branching still outperforms the learned branching rule in closed optimality gap. The experimental results seem to indicate that the former makes smarter decisions, and is on average faster in making them.
It is interesting to note that the method of \citet{Alvarez2017} performs better on homogeneous instance collections, i.e., when the problems in the training and the test set are of the same type.

\subsubsection{Online learning to branch}
The results presented in \citet{Alvarez2017} call for several reflections. First, experiments seem to indicate that performing strong branching at the top of the tree, where branching decisions have the most impact, is highly advantageous. Second, the authors already point towards more adaptation to the problem structure as a promising direction of improvement.

These ideas are studied in \citet{Khalil2016} and the follow-up work of \citet{Alvarez2016}. Independently, these authors proposed to use online learning with a strong branching initialization at the top of the tree. One key difference is that, while \citet{Alvarez2016} continue to frame learning to branch as a prediction task (i.e., predicting the strong branching score of a variable), \citet{Khalil2016} formulate the task as that of \emph{ranking}. Indeed, one does not necessarily need to accurately predict the variables' scores, but rather which variables have relatively better ones. The latter task is easier from a learning perspective.

Framing an online learning task allows to learn a specialized ML model per instance, in this case a simple linear function.  \citet{Khalil2016} analyse their learned models by studying the weight assigned to each of the features. The first question they ask is: are the learned models obtained for each of the instances similar? Their analysis concludes that there is only a weak correlation among the models. This supports the idea that adaptation plays a key role. In spite of the learned models being quite different, they were able to identify some features that are often given high importance (large absolute weight), such as the product of the pseudocosts, and data related to the constraint matrix. For a more detailed discussion of this analysis see \citet{Khalil2019}.

\subsubsection{Offline learning with structure specialization}
\citet{Alvarez2016} and \citet{Khalil2016} use a linear mapping from variable features to output. This offers the advantage of interpretability and low computational cost. The more complex GNN model proposed by \citet{Gasse2019} (see Section \ref{sec:representations} for a description) has proven to be remarkably effective in representing MILPs for the task of variable selection and beyond. The GNN architecture that we presented in Section \ref{subsec:nn} consists of a number of parametric function compositions that enable learning more complex relations between input and output. They also require a larger amount of training data, which rules out the online learning methodology previously discussed. In order to still ensure some level of specialization, \citet{Gasse2019} propose a middle point: they argue that in many realistic cases instances of the same \emph{class} are routinely solved. This justifies learning a branching rule per instance type, as a sensible trade-off between a completely general rule (such as the one in \citet{Alvarez2017}) and a completely instance-tailored rule (like in \citet{Alvarez2016, Khalil2016}). 

\citet{Gasse2019} propose training this GNN to imitate strong branching via behavioral cloning \cite{Pomerleau1991}. In short, this means that we again disregard the actual variable scores and focus on learning relative magnitudes among them. Through this approach the authors were able to outperform reliability branching, marking a breakthrough in the learning to branch literature. Building on this work, \citet{Gupta2020,Gupta2022} propose modifications of the original loss function that further improve the performance of the learned model. \citet{Seyfi2023} additionally propose a mechanism that incorporates information about past branchings into the scoring system.

\subsubsection{Towards a general branching rule}
\label{subsubsec:treeBranching}
We have discussed methodologies that specialize to certain combinatorial structures (such as \citet{Gasse2019}) or that yield a custom strategy per instance (e.g., \citet{Khalil2016}). As discussed at the beginning of this section, these methodologies are a response to the great difficulty of learning one unique policy on an heterogeneous set of instances. \citet{Zarpellon2021} argue that, while learning such a general policy poses a big challenge, it is possible to overcome it by using a representation of the search tree to inform variable selection. Their hypothesis is that there is a higher-order shared structure among MILPs, even among those with different combinatorial structure, and that this shared structure can be captured in the space of B\&B trees. To test this they define a set of features describing the state of the search which, together with variable descriptors, are mapped into scores. We discuss these features in more detail in Section~\ref{sec:representations}. The experiments in \citet{Zarpellon2021} show that, while the model of \citet{Gasse2019} struggles to learn over a heterogeneous data distribution, adding the tree context is beneficial to the generalisation performance of the model. 
\citet{Lin2022} take this idea one step further and propose to keep a record of the 
features in \cite{Zarpellon2021}  at each node where branching was performed. At every step, branching decisions are informed by this historical data, which is carefully aggregated and combined with the descriptions of the variables. 

The work in \cite{Zarpellon2021,Lin2022} highlights the potential of using solver statistics to influence branching decisions. It is unclear how the features proposed in \citet{Zarpellon2021} are used to score branching candidates,\footnote{This is because neural networks lack explanability.} but certainly this information opens the door to branching rules that switch among different behaviors at different parts of the tree, or stages of the solving process.  It is important to note that reliability branching still outperforms both \citet{Zarpellon2021} and \citet{Lin2022}, perhaps because of the overhead of computing and processing these comprehensive descriptors. Nonetheless, this work calls for further research on exploiting tree information.

\subsubsection{Expert-free learning to branch}
So far, we have discussed methods for building fast approximations of strong branching. The general consensus is that strong branching yields relatively small B\&B trees compared to other classic branching strategies, and it is therefore advantageous to try to imitate it. This idea can be challenged with two arguments. First, as pointed out by \citet{Gamrath2018}, standard implementations of strong branching benefit from using data obtained as a by-product of the score calculation, such as bound tightenings and other statistics. A rule that \emph{imitates} strong branching cannot profit from such side-effects, which means that, even in the case of perfect emulation, the imitator's performance would be worse than expected. Second, strong branching relies on the LP relaxation for scoring variables, which can provide little information in cases when the optimal LP objective value does not change with branching. In such cases strong branching cannot be considered a reliable expert. In the absence of a better alternative to strong branching, the following question arises: can we learn branching rules without expert knowledge? This question is addressed by \citet{Etheve2020} and \citet{Scavuzzo2022}, who use RL to learn branching rules from scratch.  
The computational study in \citet{Scavuzzo2022} compares one such RL methodology against the imitation learning method of \citet{Gasse2019}. 
The experiments show that on instances where strong branching performs very well, the imitation method of \citet{Gasse2019} is superior. Yet, when strong branching struggles, the RL-based method is able to find a better branching strategy, proving the point that expert-free learning is of great interest in certain cases.

The branching strategies proposed in \citet{Etheve2020} and \citet{Scavuzzo2022} do not follow demonstrations by another rule. One may therefore wonder if there is an interpretation for the learned behavior. It is unfortunately difficult to answer this question, because the architectures they use are not interpretable.

\subsubsection{Future outlook}

\paragraph{Adaptiveness.} 
The computational studies suggest that no single branching rule outperforms others universally across all instances. Consequently, a desirable approach involves a rule that dynamically adapts its behavior to the specific characteristics of each situation. One way in which one can introduce such adaptiveness is by controlling the distribution of data samples that the model uses for learning. Some authors propose to specialize to each given instance. This means that a set of parameters $\vtheta$ is generated for each new instance, obtained by allowing the ML model to only see data coming from that instance. Another possibility is to use data samples coming from instances of the same problem class. This gives us parameters $\vtheta$ that specialize to a given class and work well, on average, on different instances with shared combinatorial structure. We can therefore achieve adaptiveness on different scales. We can also understand adaptiveness of a branching rule as some sort of mechanism to change its behavior on different parts of the search process. Some progress has also been made in this regard by investigating different statistical measures that can inform branching (see Section~\ref{subsubsec:treeBranching}). Yet, little is understood about how these metrics are used or in which ways we can further enhance performance without sacrificing speed.

\paragraph{Expert guidance.} The strong branching heuristic has been used by many as an expert from which we can learn effective decision-making. The claim that strong branching is a desirable strategy to follow has recently been challenged, with some notable examples of instances where strong branching scores provide no useful information. Interesting research directions include identifying new experts, new strategies to better imitate them, or, conversely, more efficient approaches to learning without expert knowledge. The latter case is what is referred to in \citet{Bengio2021} as \emph{experience}: there is no clear mathematical understanding of what should be statistically learned (the expert), so exploration should be performed. In turn, this clearly calls for RL methods that are also natural candidates to extend the work in \citet{Zarpellon2021,Lin2022} and exploit tree information.

\paragraph{New directions.} The work we surveyed showcases the potential in mixing the extensive body of domain knowledge in variable selection with new learning techniques. Still a lot of open questions remain. For example, little attention has been directed towards highlighting important subsets of variables, as opposed to choosing a single one at each node. \citet{Khalil2022finding} propose an approach to finding such important subsets, in this case the so-called \emph{backdoors} \cite{Dilkina2009}, and show promise in using them as prioritized branching candidates. Another relevant gap is the absence of expert knowledge for certain classes of MILPs. In any case, it is clear that new ML-based methods need to build upon the pre-existing knowledge on variable selection to achieve a fruitful combination.

\subsection{Cutting planes}
\label{subsec:cutting}
Cutting plane routines are another essential part of modern MILP solvers.
They tend to work in \emph{rounds}, also called \emph{separation rounds}. Given an LP-relaxation solution $\vx^{LP}$, one round consists of generating a number of valid cuts from different families, selecting a subset of them via a selection criterion, adding them and finally resolving the LP relaxation, where $\vx^{LP}$ will now be infeasible. 
A good selection criterion is critical to improving the LP relaxation while avoiding an excessive number of cuts, which would slow down LP solving as well as lead to numerical instability. Several metrics have been proposed for the purpose of scoring cuts (see, e.g.,\@ \citet{Wesselmann2012} for an overview), 
and more recently the question of cut selection has been addressed with ML-driven predictions, which is the topic of this section. For a more in-depth discussion of ML for cut selection in MILP and beyond we refer to \citet{Deza2023}. 

As noted in Section \ref{sec:solvers}, cuts are typically more heavily applied in the root node, and for this reason the work that we survey focuses on cut selection in the root node. Still, there is no obstacle to applying these methods in other nodes. However, it is unclear whether using cuts outside the root node is computationally beneficial (see \citet{Berthold2022learning} for a discussion of this topic).

\subsubsection{Single-cut selection}
\citet{Tang2020} and \citet{Paulus2022} frame the task of cut selection as an MDP. At each step $k$, a single cut $c_k$ is selected from a cutpool $\mathcal{C}_k$, after which the LP relaxation is resolved. In particular, let $\mathcal{C}$ be a collection of cuts, and let $z(\mathcal{C})$ be the result of solving the root LP relaxation after adding all cuts in $\mathcal{C}$. The metric of interest to these authors is the 
LP-bound improvement attained by a cut $c$ at step $k$, defined as 
$$\Delta_k(c):=z\left(\{c_1,...,c_{k-1},c\} \right) - z\left(\{c_1,...,c_{k-1}\} \right) \, .$$
This MDP model is summarized in Figure \ref{fig:cutting_models}a, where we use the abbreviation $z_k=z\left(\{c_1,...,c_{k}\} \right)$.

\citet{Paulus2022} use imitation learning. Their expert is the result of explicitly calculating $\Delta_k(c)$ for each possible cut $c$ in the cutpool $\mathcal{C}_k$, and then picking the cut with highest $\Delta_k$. This means that their approach is greedy: they look at bound improvement one step ahead. Their computational studies show that this greedy heuristic that they aim to imitate is in fact very effective in improving the LP bound after $T$ cuts have been added, compared to other selection heuristics.

Another approach is that of \citet{Tang2020}, who use RL with reward $R_k=\Delta_k(c_k)$. 
Due to to the discount factor (see Equation~\ref{eq:rl}), this RL strategy offers the possibility to learn less greedy policies, and doing so without explicitly computing $\Delta_k(c)$ for each $k$ and $c\in\mathcal{C}_k$. However, many RL algorithms are known to suffer from sample inefficiency and lack of generalization \cite{Dulac2021, Kirk2023}.
\citet{Paulus2022} compare their approach to the one of \citet{Tang2020} and to SCIP's v.7.0.2 default rule. They use the following metric (lower value is better)
\begin{equation}
\label{eq:inverseIGCintegral}
    \sum_{k=1}^{T} \frac{z^* - z_k}{z^*-z_0}\,,
\end{equation}
with $T=30$ being the total number of cuts added, and $z^*$ being the pre-computed optimal solution. Notice that this metric is constructed in a way that might favor greedy policies. The computational results favor the method of \citet{Paulus2022}. The authors also show promising results in solving time when their method is incorporated into SCIP and instances are solved to optimality, though they do not include root node processing time.

\subsubsection{Multi-cut selection}
A potential criticism to the cut selection models is that solvers usually add more than one cut per round, in order to reduce the number of times the LP needs to be resolved. In fact, historically this proved crucial to the efficiency of cutting plane routines \cite{Balas1996, Cornuejols2007}. Having information about the LP solution after the addition of each cut is therefore unrealistic. \citet{Paulus2022} do not include root node processing times in their report, a metric under which their method is likely unfavored. Furthermore, metrics like the one in Eq. \ref{eq:inverseIGCintegral} encourage greedy bound improvements, whereas in practice cuts can work together to achieve a better bound improvement at the end of the round. In other words, the optimality gap closed after a full separation round is likely a more informative metric.

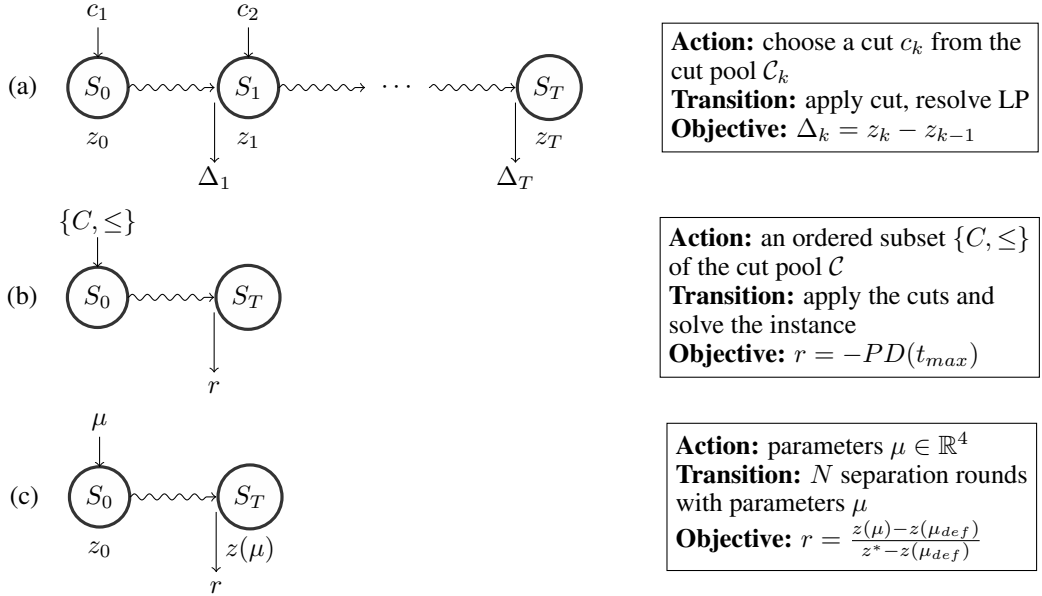
\begin{figure}[t]
\centering
    \begin{tikzpicture}
    \tikzset{decoration={snake,amplitude=.4mm,segment length=2mm, post length=.8mm,pre length=0mm}}
    \tikzstyle{state}=[{circle, draw=black!80, fill=white!95!black!5, very thick, minimum size=7mm}]

    \node at (-1,0) [](){(a)};
    \node at (0,0) [state](0){$S_0$};
    \node at (2,0) [state](1){$S_1$};
    \node at (4,0) [](dots){$\cdots$};
    \node at (6,0) [state](T){$S_T$};

    \draw[->, decorate] (0.4,0) -- (1.55,0);
    \draw[->, decorate] (2.4,0) -- (3.55,0);
    \draw[->, decorate] (4.4,0) -- (5.55,0);

    \node at (0,1)[]{$c_1$};
    \draw[->] (0,0.8) -- (0,0.4);
    \node at (2,1)[]{$c_2$};
    \draw[->] (2,0.8) -- (2,0.4);

    \node at (1.55,-1.2)[]{$\Delta_1$};
    \draw[->] (1.55,-0.2) -- (1.55,-1);
    \node at (5.55,-1.2)[]{$\Delta_T$};
    \draw[->] (5.55,-0.2) -- (5.55,-1);


    \node at (0,-0.7)[]{$z_0$};
    \node at (2,-0.7)[]{$z_1$};
    \node at (6,-0.7)[]{$z_T$};

    \node at (10,0) [draw, align=left]{
        \textbf{Action:} choose a cut $c_k$ from the \\
        cut pool $\mathcal{C}_k$ \\
        \textbf{Transition:} apply cut, resolve LP \\
        \textbf{Objective:} $\Delta_k=z_k-z_{k-1}$ 
    };
    
    \end{tikzpicture}

    \begin{tikzpicture}
    \tikzset{decoration={snake,amplitude=.4mm,segment length=2mm, post length=.8mm,pre length=0mm}}
    \tikzstyle{state}=[{circle, draw=black!80, fill=white!95!black!5, very thick, minimum size=7mm}]

    \node at (-1,0) [](){(b)};
    \node at (0,0) [state](0){$S_0$};
    \node at (2,0) [state](1){$S_T$};

    \draw[->, decorate] (0.4,0) -- (1.55,0);

    \node at (0,1)[]{$\{C,\leq\}$};
    \draw[->] (0,0.8) -- (0,0.4);

    \node at (1.55,-1.2)[]{$r$};
    \draw[->] (1.55,-0.2) -- (1.55,-1);

    \node at (10,0) [draw, align=left]{
        \textbf{Action:} an ordered subset $\{C, \leq\}$\\
        of the cut pool $\mathcal{C}$ \\
        \textbf{Transition:} apply the cuts and \\
        solve the instance \\
        \textbf{Objective:} $r=-PD(t_{max})$ 
    };
    
    \end{tikzpicture}

    \begin{tikzpicture}
    \tikzset{decoration={snake,amplitude=.4mm,segment length=2mm, post length=.8mm,pre length=0mm}}
    \tikzstyle{state}=[{circle, draw=black!80, fill=white!95!black!5, very thick, minimum size=7mm}]

    \node at (-1,0) [](){(c)};
    \node at (0,0) [state](0){$S_0$};
    \node at (2,0) [state](1){$S_T$};

    \draw[->, decorate] (0.4,0) -- (1.55,0);

    \node at (0,1)[]{$\mu$};
    \draw[->] (0,0.8) -- (0,0.4);

    \node at (1.55,-1.2)[]{$r$};
    \draw[->] (1.55,-0.2) -- (1.55,-1);


    \node at (0,-0.7)[]{$z_0$};
    \node at (2,-0.7)[]{$z(\mu)$};

    \node at (10,0) [draw, align=left]{
        \textbf{Action:}  parameters $\mu\in\mathbb{R}^4$ \\
        \textbf{Transition:} $N$ separation rounds \\
        with parameters $\mu$ \\
        \textbf{Objective:} $r=\frac{z(\mu) - z(\mu_{def})}{z^*- z(\mu_{def})}$ 
        
    };
    
    \end{tikzpicture}

\caption{Three models for learning to cut: (a) \citet{Tang2020, Paulus2022}, (b) \citet{Wang2023}, (c) \citet{Turner2023}. Here $PD(t_{max})$ refers to the primal-dual integral (see Section \ref{subsec:metrics}). }
\label{fig:cutting_models}
\end{figure}

To address the potential interactions among cuts, \citet{Wang2023} propose a policy that decides the fraction of cuts from the pool to be selected, and scores ordered subsets of this size. See Figure \ref{fig:cutting_models}b for a summary of this selection model. The authors train the policy with an RL algorithm and use end-of-run statistics like solving time as the reward.
This requires solving an MILP to optimality for each training sample. While this allows to learn a mapping from cut selection to actual solver performance (instead of the root node bound improvement, which is just a proxy for performance) this sample collection comes at a great computational cost.

A third model to learning cut selection is proposed in the work of \citet{Turner2023} (see Fig \ref{fig:cutting_models}c). The procedure builds upon SCIP's default strategy, a rule that has been carefully curated through computational studies \cite{Achterberg2007, Wesselmann2012}. This rule combines four cut scoring functions $s_i:\mathbb{R}^{n+1} \to \mathbb{R}_{+},\ i=1,\dots,4$ via a convex combination 
$$s(c, \mu)=\mu_1s_1(c)+\mu_2s_2(c)+\mu_3s_3(c)+\mu_4s_4(c)\, ,$$
$$\sum_{i=1}^4\mu_i=1, \hspace{3mm} \mu_i\geq 0, \hspace{3mm} i\in\{1,2,3,4\} \, .$$

We refer to \citet{Turner2023} for the definition of the scoring functions $s_i$, $i=1,...,4$. The problem of choosing a good set of parameters $\mu\in\mathbb{R}^4$ has been studied from a learning theoretical perspective by \citet{Balcan2021}. As a next step, \citet{Turner2023} argue, through both theoretical and computational arguments, in favor of adapting the coefficients $\mu$ to the instance being solved, as opposed to a unique, static choice. One of their experiments consists of finding a custom set of parameters $\mu$ per instance through grid search. While impractical, this experiment uncovers the large potential for improvement when adapting the value of $\mu$. In order to exploit this potential in a more realistic way, \citet{Turner2023} devise a policy that, given an instance, chooses custom parameters $\mu$, and they train it via RL. The computational study shows that the learned policy is competent in its task, outperforming random selection in terms of closed optimality gap at the root node. However, and perhaps surprisingly, their policy does not perform consistently better than SCIP's default settings when it comes to final solving time. 

\subsubsection{Beyond scoring} 
Cut scoring for selection is an essential part of cut management. Yet, there are other important decisions. \citet{Wang2023} incorporate the number of cuts added as a decision that the policy must make. Very recently, \citet{Li2023} defined a learning task for separator configuration. The objective is to select a subset of the available separators. Only the selected separators will then be active and contribute to the cutpool, meaning that this selection step happens \emph{before} the cut selection phase. \citet{Li2023} propose a methodology to overcome the high dimensionality of the configurations space, which is $2^{M\cdot R}$, with $M$ the number of separators and $R$ the number of cutting rounds. Their experimental findings show a lot of promise. More research into adapting other parametric choices could provide further insights.

\subsubsection{Future outlook}
\paragraph{Measuring performance.} What is the purpose of the cutting routines? One is inclined to believe that cuts should strengthen the 
LP relaxation, hence bringing the LP bound closer to the optimal value. However, will this result in a faster solve? \citet{Turner2023} experimentally measure the (kind of folklore) fact that a better root LP bound does not always translate into a shorter solving time. Other than the clear influence of the number of added cuts, many solver components can be affected by the cut choice resulting in different performance. \citet{Wang2023} address this by incorporating solving time as a reward signal, instead of root LP bound. However, observing final performance comes at great computational cost, which could be prohibitive for larger instances. More research is needed on how to efficiently navigate this trade-off.
\paragraph{Multi-cut rounds.} The selection model of \citet{Tang2020} and \citet{Paulus2022} adds one cut at a time, resolving the LP at each iteration. Their work constitutes an important step towards learning cut selection rules. However, in a practical setting, such procedure is unlikely to outperform models that do limited LP resolving by adding multiple cuts at once. Going forward, models like the one of \citet{Wang2023} or \citet{Turner2023} have more potential for improvement.

\subsection{Node selection}
\label{subsec:nodesel}
Primal heuristics have the clear goal of improving the best known MILP-objective value (primal bound). Analogously, branching rules and cutting routines are typically designed to improve the dual bound. Node selection policies have the difficult task of balancing both goals, which 
are often at odds. As usual, a better node selection rule is one that results in the shortest solving time. This is typically associated with smaller search trees. For that, one needs to avoid processing nodes that could be pruned if the optimal solution was known in advance. Finding a good (or even optimal) solution fast makes that task easier to accomplish. 

One strategy is to first process nodes with the best (lowest) known lower bound. This is called best first search (BFS) and has the benefit of quickly improving the dual bound. Note that for all B\&B trees there is a node selection policy of BFS type\footnote{There can be more than one BFS policy because of ties.} that minimizes the number of processed nodes \cite{Achterberg2007}.  However, nodes with good lower bounds are usually at the top of the tree whereas feasible solutions 
are typically found in deeper nodes or leaves. The depth first search (DFS) strategy prioritizes children or siblings of the node that was last processed, aiming at quickly finding feasible solutions. This has the added benefit of faster node processing times, on account of the similarity between sub-problems that are solved consecutively, which usually differ in one variable bound change. In practice, node selection rules alternate between both behaviors, while considering other statistics about branching that allow for estimating the cost of integrality.

\citet{He2014} propose to learn from a node selection oracle that always chooses to process the node that is on a path to the optimal solution. This requires knowledge of the optimal solution during the training process, but not at test time. Similarly, \citet{Yilmaz2021} prioritize nodes that contain high-quality solutions using a policy that always picks a child of the current node and uses an ML-based prediction to choose among these children. Once the dive is finished, they propose different ways in which the next node can be selected. \citet{Labassi2022} learn a function that compares any two nodes in the tree. This can be used to substitute lower bounds as the proxy of a node's potential. The learned comparison function can be used in combination with different selection strategies, such as picking the child node with highest potential. Finally, \citet{Khalil2022mipgnn} guide the search based on a prediction of the optimal solution. They look at the fixed variables at each node and measure the similarity between the fixed and predicted values. This favors nodes where the partial solution aligns with the predicted solution. 

\subsubsection{Future outlook}

The experimental results of the papers surveyed above show promise, yet the margins of improvement remain small. Often, an effective heuristic schedule and branching strategy are much more crucial and, when chosen correctly, make the impact of the node selection strategy relatively small. Still, it is interesting to observe that ML has opened new opportunities in an area where research has been pretty much inexistent for decades. This suggests that there is some potential for looking at an ``old" problem from a different perspective and with new tools. For example, one could pair different node selection strategies with restarts \citep{Anderson2019}, i.e., changing the node selection in a more dramatic way over time.

\subsection{Configuration decisions}
MILP solvers are highly parametric. To illustrate this, consider SCIP version 8.0.0: it has more than 2000 parameters that the user can tune.  A good parameter configuration that suits the instances being solved (for example a certain class of instances) can have a crucial effect on the solving process. Again, we can look at this problem with a Machine Learning lens: we can base certain parametric decisions on a prediction given by an ML model. One can see this as falling under the realm of Automatic Algorithm Configuration (AAC). However, AAC methods typically entail configuring a large number of parameters at the same time (see, e.g.\@ \citet{Hutter2009}). Here, we are interested in the use of ML to answer a single parametric question or, at least, one question at a time. Furthermore, and contrary to other AAC methodologies, the ML models make use of a description of the instance. These models are trained on a heterogeneous collection of instances and allow for instance-dependent parameter prescriptions, as opposed to a single configuration for the given instance distribution.
Some work in this area is summarized in Table~\ref{tab:other}. \\

\begin{table}[ht]
\centering
\caption{Four examples of using ML to set a solver parameter. The ML model is responsible for answering a single parametric question by choosing one of the available options. 
}
\label{tab:other}
\begin{tabular}{rccc}
  & Component & Question & Options \\ 
 \hline 
\multirow{2}{*}{\citet{Kruber2017}} & \multirow{2}{*}{General} & Should the Dantzig-Wolfe & \multirow{2}{*}{yes / no} \\
 &  & decomposition be used? &  \\
\multirow{2}{*}{\citet{Hendel2019}}  & \multirow{2}{*}{LP solver} & Which simplex pricing  & devex / steepest / \\
 &  & rule to use? & quick-start steepest \\
\multirow{2}{*}{\citet{Berthold2021}} & \multirow{2}{*}{Presolve} & Which scaling method & Standard / \\
 &  & to apply? & Curtis-Reid \\
\multirow{2}{*}{\citet{Berthold2022learning}} & \multirow{2}{*}{Cutting} & Should cuts be applied & \multirow{2}{*}{yes / no} \\
 &  & outside the root node? &  \\
 \hline
\end{tabular}
\end{table}

This avenue of research has already fostered considerable success. Notably, the method presented in \citet{Berthold2021} is used by default in Fico XPRESS v.8.9. The  great potential of customised configurations is highlighted in problems where the preprocessing techniques have a broader impact, such as in mixed integer nonlinear programming. A prime example of this is \citet{Bonami2018,bonami2022classifier}, where the authors prescribe for each mixed integer quadratic programming instance if the quadratic objective function should be linearised or not. However, note that preprocessing has been shown to be the single most impactful component of MILP solvers \cite{Achterberg2013}, hence the use of ML to configure the MILP algorithmic decisions based on the characteristics of an instance or a class of instances has a strong potential. 

\subsubsection{Future outlook}

A related line of research is that of using ML to predict search completion (see, e.g., \citet{Fischetti2019, Hendel2022b}). One can envision that in the future these predictions can be used to trigger a restart \citep{Anderson2019}, a technique that has gained a lot of attention and that allows to \emph{re}configure algorithmic decisions based on the evolution of the B\&B. In other words, ML models can be used not just during preprocessing, but also to prescribe a change in strategy during the solve, especially because before a restart some (sometimes) extensive data collection has happened, data to be leveraged by ML. Along the same lines, it is worth mentioning the recent attempts to leverage data to better solve sequences of MILP instances that differ very little one from another. This has been the focus of the 2023 MIP challenge,\footnote{\url{https://github.com/ambros-gleixner/MIPcc23}} and again it pertains to effectively configure an MILP solver, i.e., its algorithmic decisions, by exploiting data associated with instance similarities and data collected from previous runs. For the challenge, several classes of instances were proposed, where the instances in each class differ very little, for example, only in the coefficients of the objective function. The solver that won the competition \citep{patel2023progressively} was able to leverage the data of the (previous) runs, for example, the pseudocosts for making better branching decisions.

Many open questions in this vast research area still exist, making this a promising area of future work.

\section{Problem representation}
\label{sec:representations}
The standard form of Mixed Integer Linear Programs is the one presented in Equation~\ref{eq:MILP}. Given $A\in \mathbb{Q}^{m\times n}$, $\vc\in\mathbb{Q}^n$, $\vb\in\mathbb{Q}^m$ and the partition $(\mathcal{A}, \mathcal{B}, \mathcal{C})$ of the variables, an MILP solver can start solving the instance. The ML contributions to the solving process surveyed in Section~\ref{sec:tasks} also require information about the problem, but the data \{$A$, $\vc$, $\vb$, $(\mathcal{A}, \mathcal{B}, \mathcal{C})$\} may be insufficient or unfit for the task at hand. In this section we review different methodologies to construct a representation $X$ of the problem being solved. This representation is the input to a parameterized policy $f(X,\theta)$ to be trained for a specific task.

Let us start by listing the desirable properties of a representation. 
\begin{enumerate}
    \item Permutation invariance:\footnote{Note that permutation invariance is an issue beyond the ML context. The performance of MILP solvers can be affected by a change in the order of the variables or constraints. Such seemingly irrelevant changes that have an impact on the solution process are a known issue called \emph{performance variability} (see \citet{Lodi2013}).} permuting the order of the variables and/or constraints should leave the representation unchanged.
    \item Scale invariance: scale invariance is preferred to keep values within controlled ranges, which helps the learning process. This can be achieved with a normalization step.
    \item Size invariance: the size of the representation should not depend on the size of the instance. This is, we require a fixed-sized description of each element that needs representation, e.g., each variable or each node.
    \item Low computational cost: low cost of extracting, storing and processing the data.
\end{enumerate}

In the following, we will make a distinction with respect to descriptors that represent general properties of the MILP and descriptors that relate to a specific variable or constraint. These descriptors may be static in nature or they may dynamically change during the solving process. We will also discuss global descriptions of the process versus local (subproblem) properties. Note that some approaches use no description of the instance, and instead learn exclusively from the performance metric (see, e.g., \citet{Chmiela2021, Hendel2022}).

\subsection{Representing variables individually}
Some of the learning tasks discussed in Section~\ref{sec:tasks} require a description of each variable individually. Clearly, branching rules fall under this category. This is also the case for prediction-driven heuristics and neighborhood selection policies, for which variables are mapped to values or probabilities. In this section, we will discuss three important approaches to building variable representations and how they relate to the different approaches of Section~\ref{sec:tasks}. 

A straightforward approach to building variable representations is to gather a number of descriptors into a simple \emph{vector representation} for each variable. \citet{Khalil2016} propose a number of such descriptors, including different statistics about the set of constraints in which each variable participates. These statistics aggregate information whose length would otherwise depend on the problem size. For example, for a variable $j\in\mathcal{I}$, using the constraint coefficients $\{a_{ij}\}_{i=1}^m$ would yield a vector whose length depends on $m$, which is undesirable. On the contrary, using the average of these coefficients gives a size-independent descriptor. This is a necessary step but calls to question which statistics should be included or excluded in this feature engineering step. 

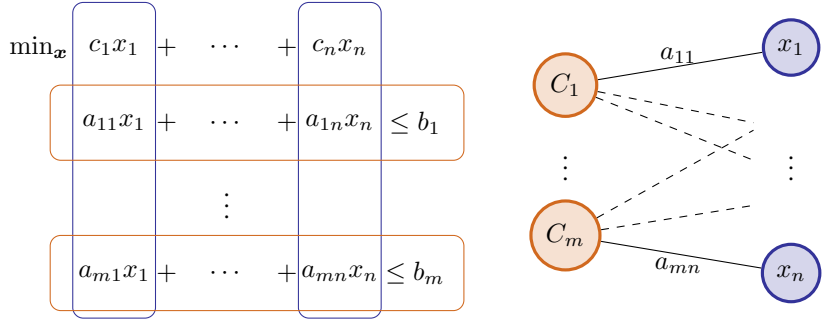
\begin{figure}[t]
\centering
    \begin{tikzpicture}
    \tikzset{decoration={snake,amplitude=.4mm,segment length=2mm, post length=.8mm,pre length=0mm}}
    \tikzstyle{varnode}=[{circle, draw=prettyblue, fill=prettyblue!20, very thick, minimum size=7mm}]
    \tikzstyle{consnode}=[{circle, draw=cocoabrown, fill=cocoabrown!20, very thick, minimum size=7mm}]

    \node at (0,0) [](){$\min_{\vx}$};
    \node at (1,0) [](){$c_1x_1$};
    \node at (1.7,0) [](){$+$};
    \node at (2.5,0) [](){$\cdots$};
    \node at (3.3,0) [](){$+$};
    \node at (4,0) [](){$c_nx_n$};

    \node at (1,-1) [](){$a_{11}x_1$};
    \node at (1.7,-1) [](){$+$};
    \node at (2.5,-1) [](){$\cdots$};
    \node at (3.3,-1) [](){$+$};
    \node at (4,-1) [](){$a_{1n}x_n$};
    \node at (5,-1) [](){$\leq b_1$};

    \node at (2.5,-2) [](){$\vdots$};

    \node at (1,-3) [](){$a_{m1}x_1$};
    \node at (1.7,-3) [](){$+$};
    \node at (2.5,-3) [](){$\cdots$};
    \node at (3.3,-3) [](){$+$};
    \node at (4,-3) [](){$a_{mn}x_n$};
    \node at (5,-3) [](){$\leq b_m$};

    \node[rounded corners,
    draw = prettyblue,
    minimum width = 1.1cm, 
    minimum height = 4.2cm] () at (1,-1.5){};
    \node[rounded corners,
    draw = prettyblue,
    minimum width = 1.1cm, 
    minimum height = 4.2cm] () at (4,-1.5){};

    \node[rounded corners,
    draw = cocoabrown,
    minimum width = 5.5cm, 
    minimum height = 1cm] () at (2.9,-1){};
    \node[rounded corners,
    draw = cocoabrown,
    minimum width = 5.5cm, 
    minimum height = 1cm] () at (2.9,-3){};

    \node at (7,-0.5) [consnode](c1){$C_1$};
    \node at (7,-1.5) [](){$\vdots$};
    \node at (7,-2.5) [consnode](cm){$C_m$};
    \node at (10,0) [varnode](x1){$x_1$};
    \node at (10,-1.5) [](){$\vdots$};
    \node at (10,-3) [varnode](xn){$x_n$};
    \draw[-] (c1) -- (x1);
    \draw[dashed] (c1) -- (9.5,-1);
    \draw[dashed] (c1) -- (9.5,-1.5);
    \draw[dashed] (cm) -- (9.5,-1.1);
    \draw[dashed] (cm) -- (9.5,-2.1);
    \draw[-] (cm) -- (xn);
    \node at (8.5,-0.1) [](){$a_{11}$};
    \node at (8.5,-2.9) [](){$a_{mn}$};
    
    \end{tikzpicture}

\caption{The bipartite graph representation of an MILP.}
\label{fig:bipartite}
\end{figure}

Alternatively, one can use the \emph{bipartite graph representation} of the MILP. This graph, shown in Figure~\ref{fig:bipartite}, 
is constructed as follows: each constraint and each variable have a corresponding representative node. A constraint node is connected to a variable node if the corresponding variable has a non-zero coefficient in the corresponding constraint. Each node has an associated vector descriptor. The advantage of using a graph representation is that this data structure can be parsed by a Graph Neural Network (GNN, see Section~\ref{subsec:nn} for a formal definition). This type of architecture automatically handles inputs of different sizes.  The data aggregation step that ensures size invariance is one of the learnable mappings (function \texttt{comb} in Definition \ref{def:gnn}). This is, instead of manually engineering a mechanism to aggregate information coming from the constraints, this mechanism is automatically \emph{learned}. The use of GNNs for combinatorial optimization has experienced a rise in popularity in recent years \citep{Cappart2023} because of their ability to capture the structural properties of the instances without excessive engineering. 

\citet{Gasse2019} introduced
this bipartite graph representation to make predictions about variables, using a GNN as a mapping. Their node descriptors include structural information, such as constraint coefficients, that capture how variables interact with each other through constraints. Other than this, both \citet{Khalil2016} and \citet{Gasse2019} include information about the LP solution and other basic variable features such as their objective coefficient or variable type. \citet{Zarpellon2021} take a different perspective, stressing the importance of historical data 
collected during the B\&B tree. This strategy resembles SCIP's default branching rule, which considers information about past branchings, collected conflicts or cutoffs. The representation used in \citet{Zarpellon2021} includes this variable information and, additionally, global information about the search tree. They argue that such a description can uncover shared structures among very diverse MILPs (see Section~\ref{subsec:branching}). The information collected includes statistics about the node being processed, tree composition and shape, and bound statistics, with a particular focus on unprocessed nodes. 

These three approaches to MILP variable representation are compared in Table \ref{tab:var_rep_model}. We report on features that describe variables individually, therefore excluding the complete list of tree features of \citet{Zarpellon2021} (which can be found in their appendix). Table~\ref{tab:var_rep_model} showcases that, while some common features exist, the different representations have distinct focuses.

Most of the work surveyed in Section~\ref{sec:tasks} uses the graph representation of \citet{Gasse2019}, either exactly (e.g., \cite{Sonnerat2021, Gupta2022, Scavuzzo2022}) or with small modifications of the variable descriptors (e.g., \cite{Wu2021, Liu2022, Khalil2022mipgnn}). Hybrid models also exist. In particular, \citet{Gupta2020} propose extracting the \citet{Gasse2019} representation at the root node and the \citet{Khalil2016} representation in the rest of the nodes of the search tree. The reason is that, while the graph representation is rich, it is also computationally expensive. \citet{Gasse2019} overcome this by using a GPU (graphics processing unit) to accelerate the execution of their learned function. Such computation on GPUs is common practice in the ML community, but one could argue that it is unrealistic to require the availability of a GPU. The hybrid model of \citet{Gupta2020} reuses the rich but expensive representation of the root node in combination with the features of \citet{Khalil2016} that are cheap and update the description at every node. This proves very effective, with their best performing model outperforming both reliability branching and the model of \citet{Gasse2019} when executed without GPU acceleration. 

\begin{table}[tb]
\centering
\caption{Three key approaches to variable descriptors.}
\label{tab:var_rep_model}
\begin{tabular}{rccc}
  & \citeauthor{Khalil2016} & \citeauthor{Gasse2019} & \citeauthor{Zarpellon2021} \\ 
Basic    & & &    \\ \hline \hline 
objective coefficient & \ding{52} & \ding{52} &  \\
upper/lower bound &  & \ding{52} &  \\
type &  & \ding{52} &  \\
& & & \\
Structural & & & \\ \hline \hline 
\# of constraints the variable is in  & \ding{52} & implied & \\
min and max ratios $a_{ij}$ to $b_i$ & \ding{52} & implied & \\
min and max ratios $a_{ij}$ to $\sum_k a_{kj}$ & \ding{52} & implied & \\
stats for constraint degrees & \ding{52} & implied & \\
stats for constraint coefficients &   &   & \\
stats for active constraints & \ding{52} & implied & \\
& & & \\
LP solution & & & \\ \hline \hline 
LP basis status & & \ding{52} & \\
LP sol value & & \ding{52} & \ding{52}\\
LP sol at bound & & \ding{52} & \\
LP sol fractionality & \ding{52} & \ding{52} & \\
LP sol reduced cost & & \ding{52} & \\
LP sol age & & \ding{52} & \\
Root LP sol value & & \ding{52} & \\
& & & \\
Incumbent & & & \\ \hline \hline 
incumbent value & & \ding{52} & \\
average incumbent value & & \ding{52} & \ding{52} \\
& & & \\
Tree statistics & & & \\ \hline \hline 
average branching depth & & & \ding{52} \\
conflict score & & & \ding{52} \\
conflict length score & & & \ding{52} \\
average conflict length & & & \ding{52} \\
pseudocost score & \ding{52} & & \ding{52} \\
pseudocost count stats & & & \ding{52} \\
average inference score & & & \ding{52} \\
average number of inferences  & & & \ding{52} \\
average cutoff score & & & \ding{52} \\
average cutoffs of variable & \ding{52} & & \ding{52} \\
\# of implications derived  & & & \ding{52} \\
\# of cliques of variable & & & \ding{52} \\
\end{tabular}
\end{table}

\subsection{Representing constraints individually}
Analogously to variables, a description of the problem's constraints may be needed. Here, we refer both to original problem constraints and additional valid constraints that can be added as cuts. A constraint representation is necessary in two cases. First and undoubtedly, whenever the task requires a decision over said constraints (e.g., which cuts to add). Second, the bipartite graph representation of MILPs discussed in the previous section (see Figure~\ref{fig:bipartite}) also calls for a description of the constraints, even in the case where they are to be aggregated at a later stage.

\citet{Gasse2019} first proposed the bipartite graph representation using a small number of descriptors for the constraint nodes. In particular, they use the cosine similarity\footnote{The cosine similarity between two vectors $\va$ and $\vb$  of the same length is defined as $\va \cdot \vb / ||\va||_2 ||\vb||_2$, that is, their dot product divided by the product of their norms.} with the objective coefficients, the constraint right-hand-side, and LP information such as basis status and dual bound. This type of concise descriptions of the constraints is frequently used for the learning tasks associated with branching or primal heuristics.

In the case of cut selection, a more detailed description is preferred. \citet{Paulus2022} extend the graph representation of \citet{Gasse2019} with metrics that are typically considered in cut selection, such as violation, objective parallelism or sparsity. \citet{Wang2023} describe each cut with a single vector of classical cut scores (see. e.g., \cite{Wesselmann2012}). In the case of \citet{Turner2023}, the classical cut scores are intrinsically taken into account in the definition of their learning task (see Section~\ref{subsec:cutting}). For this reason, they use a graph representation of the model with a small amount of variable and constraint descriptors. 

\subsection{Representing a (sub-)MILP}
To conclude, it is worth observing that some decision tasks are formulated at an instance or node level, therefore needing a global representation of the MILP and perhaps also the solving process. One possible approach is to aggregate variable and constraint descriptors coming from the bipartite graph representation. \citet{Liu2022} use the average of the variable descriptors, while \citet{Labassi2022} concatenate the average variable descriptors and the average constraint descriptors.

Other approaches include the one of \citet{Khalil2017}, who, in the context of scheduling of primal heuristics, build a vector representation of the current node. This representation includes comparisons to the root node and context on the node's position within the tree. Great focus is put on information coming from the LP solution, such as the objective value, average fractionalities, and statistics on the constraint activity. Berthold et al. \cite{Berthold2021, Berthold2022learning} build representations that are more specialized to the particular configuration task. We refer to their work for a more detailed discussion of the problem descriptors.

\subsection{Notes on architecture}

\paragraph{Expressivity versus speed.} In the context of machine learning, model size refers to the number of parameters and operations that define the mapping function $f(\cdot, \vtheta)$. Larger models allow us to learn more complex relationships between input and output. However, there is a clear trade-off between model size and computational cost of execution and training.  In stark contrast to other fields of application of ML, such as large language models \cite{chatgpt} or strategic game playing \cite{alphastar} where the goal is super-human performance, the computational cost per execution of the ML model is decisive to whether or not it will beat its competitors, i.e., already highly efficient optimization software.  Navigating this trade-off is an active field of research, with promising results in both small and large models, as well as in compressing ML models without compromising the accuracy of their predictions.

\paragraph{The case for Graph Neural Networks.} 
GNNs offer a powerful representation tool for MILP. They enable instance parsing with less feature engineering, as well as size and permutation invariance. Computationally, they have shown excellent performance across different tasks. Yet, we have limited understanding of the reasons behind this success. Some recent studies have uncovered some of the factors that contribute to the success of GNNs.

\citet{Chen2023} study the separation and representation power of GNNs for LP. In particular, they study this in the context of three prediction tasks: predicting feasibility, boundedness, and the optimal solution vector for LP. The separation power is the ability of GNNs to distinguish different instances, i.e., their ability to output different results when given two different instances as input. \citet{Chen2023} prove that given two LPs, if no GNN\footnote{The authors consider the family of functions defined in Equation \ref{eq:gnns}, where the combination and aggregation functions are feed-forward neural networks.} can distinguish them, then both LPs have the same status in terms of feasibility and boundedness. Furthermore, they both have the same minimum-$l_2$-norm optimal solution up to a permutation. Finally, they also show that the three tasks mentioned above can in fact be approximated using GNNs.
Continuing this line of work, \citet{Qian2023} prove that GNNs can be used to reproduce interior point methods. In particular, they show that there exists a GNN using $\mathcal{O}(m)$ message-passing operations (see Definition \ref{def:gnn}), with $m$ the number of constraints, that can replicate any one iteration of the algorithm by \citet{Nocedal2006}. They show the same result for the more practical algorithm by \citet{Gondzio2012}. Notably, this is true when representing the LP using a modified version of the aforementioned bipartite graph representation (see Figure \ref{fig:bipartite}), where a new node is added and connected to all variable and constraint nodes. This global node adds alternative routes of communication among constraint and variable nodes and is said to represent the objective function. \citet{Qian2023} also provide a computational comparison of different GNN implementations, i.e., different $\texttt{comb}$ and $\texttt{aggr}$ functions (see Eq.~\ref{eq:gnns}). This is also an active area of research, with recent works advocating for the so-called graph attention networks (e.g., \cite{Lin2022, Seyfi2023}) and other sophisticated architectures and training methods.

These results shed some light on the representation power of GNNs for MILP and strengthen the case for using them in the context of optimization problems.

\section{Datasets and software}
\label{sec:benchmarks}
\newcommand{\ExternalLink}{%
    \tikz[x=1.2ex, y=1.2ex, baseline=-0.05ex]{%
        \begin{scope}[x=1ex, y=1ex]
            \clip (-0.1,-0.1) 
                --++ (-0, 1.2) 
                --++ (0.6, 0) 
                --++ (0, -0.6) 
                --++ (0.6, 0) 
                --++ (0, -1);
            \path[draw, 
                line width = 0.5, 
                rounded corners=0.5] 
                (0,0) rectangle (1,1);
        \end{scope}
        \path[draw, line width = 0.5] (0.5, 0.5) 
            -- (1, 1);
        \path[draw, line width = 0.5] (0.6, 1) 
            -- (1, 1) -- (1, 0.6);
        }
    }
    
The modelling power of MILP makes it a suitable language for a large range of applications. 
With the goal of measuring the performance of different algorithms, the MILP research community has curated large benchmarks, such as MIPLIB \cite{Gleixner2021}, that provide a heterogeneous set of instances coming from diverse applications. It should be noted that these 
 benchmarks are considered large for MILP standards, but are orders of magnitude smaller than typical ML benchmarks.
In the light of the methods surveyed in Section~\ref{sec:tasks}, which combine an ML component with classical optimization, there is a renewed need for collections of MILP instances. In this section we provide an overview of the collections that have been used in the growing body of literature. We restrict our discussion to benchmarks that are publicly available or whose generation code is easily accessible. 

 ML methodologies usually require vast amounts of data. For this reason, it is common to resort to instance generators, which complement the existing instance collections. Tables \ref{tab:collections} and \ref{tab:generators} provide a summary of both commonly used instance collections (with their size specification) and instance generators. 

 \begin{table}[b]
\centering
\caption{A list of some of the most relevant instance collections 
}
\label{tab:collections}
\begin{tabular}{rcccc}
Benchmark & Composition & Size & Source & URL \\ \hline
MIPLIB 2017   &  Mixed & 240 & \cite{Gleixner2021} & \href{https://miplib.zib.de/}{\ExternalLink} \\
Cor@l         &  Mixed & 364 & - & \href{https://coral.ise.lehigh.edu/data-sets/mixed-integer-instances/}{\ExternalLink} \\
NN verification   &  Homogeneous & 3692 & \cite{Nair2020} & \href{https://github.com/deepmind/deepmind-research/tree/master/neural_mip_solving}{\ExternalLink} \\
ML4CO\_1           &  Homogeneous & 10,000 & \cite{Gasse2022} & \href{https://github.com/ds4dm/ml4co-competition/blob/main/DATA.md}{\ExternalLink} \\
ML4CO\_2           &  Homogeneous & 10,000 & \cite{Gasse2022} & \href{https://github.com/ds4dm/ml4co-competition/blob/main/DATA.md}{\ExternalLink}  \\
ML4CO\_3           &  Homogeneous & 118 & \cite{Gasse2022} & \href{https://github.com/ds4dm/ml4co-competition/blob/main/DATA.md}{\ExternalLink} \\
\end{tabular}
\end{table}

\begin{table}[ht]
\centering
\caption{A list of some of the most relevant instance generators}
\label{tab:generators}
\begin{tabular}{ccc}
Benchmark & Problem type(s) & Source \\ \hline
\multirow{4}{*}{Tang et al.}         &  Max-cut, & \multirow{4}{*}{\cite{Tang2020}}  \\
 & Planning, & \\
 & Packing, & \\
 & Bin packing & \\ \hline
\multirow{4}{*}{Ecole}         &  Set cover, & \multirow{4}{*}{\cite{Prouvost2020}}  \\
 & Combinatorial auctions, & \\
 & Maximum independent set, & \\
 & Capacitated facility location& \\ \hline
\multirow{9}{*}{MIPLearn}         &  Bin packing, & \multirow{9}{*}{\cite{Santos2023}} \\
& Multi-dimensional knapsack, & \\
& Capacitated p-median, & \\
& Set cover, & \\
& Set packing, & \\
& Stable set, & \\
& Traveling salesman, & \\
& Unit commitment, & \\
& Vertex cover & \\ \hline
GISP         & Generalized independent set & \cite{Colombi2017}  \\ \hline
FCMNF        & Capacitated fixed-charge network flow & \cite{Hewitt2010} \\ \hline
\end{tabular}
\end{table}

Apart from the size, there is the consideration of the composition of instances. When implementing a learning-augmented solver component a specification needs to be made regarding the instances of interest. In machine learning terms, we typically talk of an \emph{instance distribution}, where the instances that are outside the scope of interest are assigned a zero probability of occurring. For some applications, it can be assumed that the representative instances have a shared combinatorial structure. The machine learning model is then expected to specialize to this structure. Conversely, some approaches are designed to detect patterns across instances of any class. Throughout Section \ref{sec:tasks} we have surveyed examples of both situations. Table \ref{tab:benchmarks_mixed_vs_hom} summarizes which of the discussed approaches uses a collection of instances with mixed structures (mixed), and which use the assumption that all instances belong to the same class (homogeneous). From this we can observe that configuration decisions are more naturally framed over mixed instance collections than other tasks, like the more complex matter of branching, where some instance specialization seems valuable. Table \ref{tab:benchmarks_hom} further shows an overview of the homogeneous datasets used in the work presented in Section \ref{sec:tasks}. We can observe a pattern that highlights differences in the instances, based on which task is more challenging. Instances like GISP or FCMNF are more commonly chosen as a challenging test bed for primal heuristics, indicating that for these problems the difficulty lies in finding (optimal) solutions. On the contrary, branching is usually tested on instances where proving optimality is the main challenge, such as the ones provided by Ecole.

\begin{table}[ht]
\centering
\caption{Classification of the literature based on whether they use a mixed or an homogeneous instance collection (or both).}
\label{tab:benchmarks_mixed_vs_hom}
\begin{tabular}{rcc}
                     & Mixed & Homogeneous \\ \toprule
\citet{Ding2020}     & & \checkmark  \\
\citet{Nair2020}     & \checkmark & \checkmark  \\
\citet{Khalil2022mipgnn}   & & \checkmark \\ \hline
\citet{Song2020}     & & \checkmark  \\
\citet{Wu2021}       & & \checkmark  \\
\citet{Sonnerat2021} & \checkmark & \checkmark \\
\citet{Liu2022}      & \checkmark & \checkmark  \\
\citet{Huang2023}   & & \checkmark \\ \hline
\citet{Khalil2017}   & \checkmark & \checkmark  \\
\citet{Chmiela2021}  & & \checkmark  \\
\citet{Hendel2022}   & \checkmark &  \\  \bottomrule

\multicolumn{3}{c}{Primal heuristics}        \\ 
   &  &  \\ 
                     & Mixed & Homogeneous \\ \hline
\citet{Khalil2016}   & \checkmark &  \\
\citet{Gasse2019}    & & \checkmark  \\
\citet{Gupta2020}    & & \checkmark  \\
\citet{Etheve2020}   & & \checkmark  \\
\citet{Nair2020}     & \checkmark & \checkmark \\
\citet{Zarpellon2021}  & \checkmark & \\
\citet{Gupta2022}    & & \checkmark  \\
\citet{Scavuzzo2022}  & & \checkmark  \\ \hline
\multicolumn{3}{c}{Branching} \\
   &  &  \\ 
                     & Mixed & Homogeneous \\ \hline
\citet{Tang2020}     & & \checkmark  \\
\citet{Paulus2022}   & & \checkmark  \\
\citet{Wang2023}     & \checkmark & \checkmark \\
\citet{Turner2023}   & \checkmark &  \\ 
\citet{Li2023}       & \checkmark & \checkmark  \\ \hline
\multicolumn{3}{c}{Cut selection}        \\
   &  &  \\ 
                     & Mixed & Homogeneous \\ \hline
\citet{He2014}       &   & \checkmark   \\ 
\citet{Yilmaz2021}   &   & \checkmark   \\ 
\citet{Labassi2022}  &   & \checkmark   \\ \hline
\multicolumn{3}{c}{Node selection}        \\

   &  &  \\ 
                     & Mixed & Homogeneous \\ \hline
\citet{Kruber2017}     & \checkmark &  \\
\citet{Hendel2019}   & \checkmark &   \\
\citet{Berthold2021}   & \checkmark &  \\
\citet{Berthold2022learning}   & \checkmark & \\ \hline
\multicolumn{3}{c}{Configuration decisions}        \\
\end{tabular}
\end{table}

\begin{table}[ht]
\centering
\caption{Common homogeneous instance collections and where they are used.}
\label{tab:benchmarks_hom}
\begin{tabular}{rcccccc}
                    & Ecole & NNv & GISP & FCMNF & Tang et al. & Other \\ \toprule
\citet{Ding2020}     & & & & & & \checkmark  \\
\citet{Nair2020}     & & \checkmark & & & & \checkmark \\
\citet{Khalil2022mipgnn}   & & & \checkmark & \checkmark & & \\ \hline
\citet{Song2020}     & & & & & & \checkmark  \\
\citet{Wu2021}       & & & & & & \checkmark  \\
\citet{Sonnerat2021} & & \checkmark & & & & \checkmark \\
\citet{Liu2022}      & \checkmark & & \checkmark & & & \\
\citet{Huang2023}   & \checkmark & & & & & \checkmark \\ \hline
\citet{Khalil2017}   & & & \checkmark & & &  \\
\citet{Chmiela2021}  & & & \checkmark & & &  \\ \bottomrule
\multicolumn{7}{c}{Primal heuristics}        \\ 
 & &   & &  & \\ 
                    & Ecole & NNv & GISP & FCMNF & Tang et al. & Other \\ \toprule
\citet{Gasse2019}    & \checkmark & & & & &  \\
\citet{Gupta2020}    & \checkmark & & & & &  \\
\citet{Etheve2020}   & & & & & & \checkmark  \\
\citet{Nair2020}     & & \checkmark & & & & \checkmark \\
\citet{Gupta2022}    & \checkmark & & &  & &  \\
\citet{Scavuzzo2022} & \checkmark & & &  & &  \\ \bottomrule
\multicolumn{7}{c}{Branching} \\
 & &  & & & \\
                     & Ecole & NNv & GISP & FCMNF & Tang et al. & Other \\ \toprule
\citet{Tang2020}     & & & & & \checkmark &   \\
\citet{Paulus2022}   & & \checkmark & & & \checkmark &  \\
\citet{Wang2023}     & \checkmark & & & & & \checkmark \\ 
\citet{Li2023}       & \checkmark & \checkmark & & & \checkmark & \checkmark \\\bottomrule
\multicolumn{7}{c}{Cut selection}        \\
& &  & & & \\ 
                     & Ecole & NNv & GISP & FCMNF & Tang et al. & Other \\ \toprule
\citet{He2014}       & & & & & & \checkmark \\
\citet{Yilmaz2021}   & \checkmark & & & & & \\
\citet{Labassi2022}  & & & \checkmark & \checkmark & & \checkmark   \\ \bottomrule
\multicolumn{7}{c}{Node selection}        \\
& &  & & & \\ 
\end{tabular}
\end{table}

\subsection{Software}
In connection to instance generators, there has been increasing interest in developing libraries that help the process of data generation, training and testing in the context of ML-augmented MILP solving. Some examples of these are the library Ecole \cite{Prouvost2020}, or the more recent MIPLearn \cite{Santos2023}. Their goal is to provide a standardized platform for the research community for fast prototyping and testing by removing the barrier of challenging software implementation. These libraries provide ways to easily implement learning tasks, such as branching or warm-starting. For an up-to-date specification of the features they provide we refer to their documentation.

\section{Conclusions}
\label{sec:conclusions}
The work covered in this article testifies to the growing interest in the integration of ML methodologies within MILP solvers. This is an emerging technology that has already fostered remarkable success within its short history, and is likely to play a key role in future algorithmic developments. Beyond the discussion of the literature, we have highlighted some methodological trends and characterised the common grounds with respect to instance representation, learning algorithms and benchmarking. 

Meaningful steps forward have been taken in answering the more pressing research questions. For example, the literature shows that some learning tasks seem to be formulated more naturally than others over heterogeneous instance collections. In other cases, an argument can be made in favour of the applicability of specializing to a certain structure, which makes the learning task easier. Studies like the ones in \citet{Zarpellon2021,Fischetti2019} indicate that instance representations that describe the global solution process allow to more easily recognize patters across different combinatorial structures. This is especially promising because a key challenge for the integration of ML-augmented methods into MILP solvers is their generalization properties, at least as long as solvers are conceived as one-configuration-fits-all software.

It is also interesting to note the various efforts to define efficient success metrics for the different learning tasks, such as imitation targets or reward functions. Solving instances to optimality is to be avoided because of the computational effort, but performance proxies that substitute solving time must be carefully chosen.

Already substantial progress has also been made in creating the right environment for easily implementing and testing ML models for MILP solving. Existing software infrastructure includes, for example, curated instance generators, code that simplifies the solver interface and standardized testing procedures. This further helps in evaluating and comparing the different methodologies. Other efforts to bring the research community together are competitions, like ML4CO \cite{Gasse2022}, which encourage progress in well-defined tasks as well as fair comparisons among the proposed methods. A significant challenge resides on the software versus hardware side: many learning methods, e.g., those relying on neural networks, especially benefit from the use of GPUs, while MILP technology is inherently CPU based. The CPU versus GPU interaction is currently a relevant obstacle for ML-augmented MILP.

Overall, a key trend seems to be building more dynamic solving strategies. MILP solvers generate a plethora of statistics during execution that often go unused. This is fertile ground for learning algorithms, which can unlock more dynamic solvers that automatically adapt the solving strategy based on prescriptions derived from such solving statistics. This poses exciting new questions and challenges.

\section*{Acknowledgements}
This work was partially supported by OPTIMAL, 
a project funded by the Dutch Research Council (NWO) under grant OCENW.GROOT.2019.015; and by TAILOR, a project funded by EU Horizon 2020 research and innovation programme under grant~952215.
The third author would like to warmly thank the Government of Canada for establishing the Canada Excellence Research Chair (CERC) program. It generously supported the ``Data Science for Real-time Decision-making" CERC at Polytechnique Montr\'eal that the author had the privilege to lead (2015--2022) and whose fantastic team has been instrumental to shape his knowledge of the topic.

\bibliographystyle{plainnat}     
\bibliography{refs}

\end{document}